\documentclass[latter,12pt]{article}
\usepackage{epsfig}
\usepackage{mathrsfs}
\usepackage{amssymb}
\usepackage{amsfonts}
\usepackage{amsmath}
\usepackage{bbm}
\usepackage{tikz}
\usepackage{graphicx}
\usepackage{subcaption}
\usepackage{sidecap}
\newskip\humongous \humongous=0pt plus 1000pt minus 1000pt

\newif\ifdtup

\catcode`\@=11

\@addtoreset{equation}{section}
\def\theequation{\thesection.\arabic{equation}}

\def\@normalsize{\@setsize\normalsize{15pt}\xiipt\@xiipt
\abovedisplayskip 14pt plus3pt minus3pt%
\belowdisplayskip \abovedisplayskip
\abovedisplayshortskip \z@ plus3pt%
\belowdisplayshortskip 7pt plus3.5pt minus0pt}

\def\small{\@setsize\small{13.6pt}\xipt\@xipt
\abovedisplayskip 13pt plus3pt minus3pt%
\belowdisplayskip \abovedisplayskip
\abovedisplayshortskip \z@ plus3pt%
\belowdisplayshortskip 7pt plus3.5pt minus0pt
\def\@listi{\parsep 4.5pt plus 2pt minus 1pt
     \itemsep \parsep
     \topsep 9pt plus 3pt minus 3pt}}

\relax

\catcode`@=12

\catcode`\@=11

\def\section{\@startsection{section}{1}{\z@}{3.5ex plus 1ex minus
   .2ex}{2.3ex plus .2ex}{\large\bf}}

\def\thesection{\arabic{section}}
\def\thesubsection{\arabic{section}.\arabic{subsection}}

\def\appendix{\setcounter{section}{0}
 \def\thesection{Appendix \Alph{section}}
 \def\thesubsection{\Alph{section}.\arabic{subsection}}
 \def\theequation{\Alph{section}.\arabic{equation}}}

\def\SymBoxes#1#2#3#4{\newdimen\un@t \un@t#3%
\raisebox{#1}{\rule{#2\un@t}{#4}\hskip-#2\un@t
\@tempdimb\un@t \advance\@tempdimb by-#4\@tempcntb#2\relax%
\@whilenum{\@tempcntb>0}\do{
\rule{#4}{\un@t}\hskip\@tempdimb \advance\@tempcntb by\m@ne}%
\hskip-#2\un@t \rule[\un@t]{#2\un@t}{#4}%
\rule[\un@t]{#4}{#4}\hskip-#4
\rule{#4}{\un@t}}\hskip-#4}                

\begin{document}

\newcommand{\beq}{\begin{equation}}
\newcommand{\eeq}{\end{equation}}
\newcommand{\bea}{\begin{eqnarray}}
\newcommand{\eea}{\end{eqnarray}}
\newcommand{\beas}{\begin{eqnarray*}}
\newcommand{\eeas}{\end{eqnarray*}}
\newcommand{\defi}{\stackrel{\rm def}{=}}
\newcommand{\non}{\nonumber}
\newcommand{\bquo}{\begin{quote}}
\newcommand{\enqu}{\end{quote}}
\renewcommand{\(}{\begin{equation}}
\renewcommand{\)}{\end{equation}}
\def\IZ{{\mathbb Z}}
\def\IR{{\mathbb R}}
\def\IC{{\mathbb C}}
\def\IQ{{\mathbb Q}}

\def\g{\gamma}
\def\m{\mu}
\def\n{\nu}
\def\a{\alpha}
\def\b{\beta}

\def\CM{{\mathcal{M}}}
\def\dCM{{\left \vert\mathcal{M}\right\vert}}

\def \d{\textrm{d}}
\def \p{\partial}

\def \Pf{\rm Pf\ }

\def \pr{\prime}

\def\Tr{ \hbox{\rm Tr}}
\def\half{\frac{1}{2}}

\def \eqn#1#2{\begin{equation}#2\label{#1}\end{equation}}
\def\de{\partial}
\def\Tr{ \hbox{\rm Tr}}
\def\H{ \hbox{\rm H}}
\def\HE{ \hbox{$\rm H^{even}$}}
\def\HO{ \hbox{$\rm H^{odd}$}}
\def\K{ \hbox{\rm K}}
\def\Im{ \hbox{\rm Im}}
\def\Ker{ \hbox{\rm Ker}}
\def\const{\hbox {\rm const.}}
\def\o{\over}
\def\im{\hbox{\rm Im}}
\def\re{\hbox{\rm Re}}
\def\bra{\langle}\def\ket{\rangle}
\def\Arg{\hbox {\rm Arg}}
\def\Re{\hbox {\rm Re}}
\def\Im{\hbox {\rm Im}}
\def\exo{\hbox {\rm exp}}
\def\diag{\hbox{\rm diag}}
\def\longvert{{\rule[-2mm]{0.1mm}{7mm}}\,}
\def\a{{\textsl a}}
\def\dag{{}^{\dagger}}
\def\tq{{\widetilde q}}
\def\p{{}^{\prime}}
\def\W{W}
\def\N{{\cal N}}
\def\hsp{,\hspace{.7cm}}
\newcommand{\C}{\ensuremath{\mathbb C}}
\newcommand{\Sp}{\ensuremath{\mathbb S}}
\newcommand{\Z}{\ensuremath{\mathbb Z}}
\newcommand{\R}{\ensuremath{\mathbb R}}
\newcommand{\rp}{\ensuremath{\mathbb {RP}}}
\newcommand{\cp}{\ensuremath{\mathbb {CP}}}
\newcommand{\vac}{\ensuremath{|0\rangle}}
\newcommand{\vact}{\ensuremath{|00\rangle}}
\newcommand{\oc}{\ensuremath{\overline{c}}}
\newcommand{\sgn}{\mathop{\mathrm{sgn}}}

\def\M{\mathcal{M}}
\def\F{\mathcal{F}}
\def\d{\textrm{d}}

\def\eps{\epsilon}

\begin{flushright}
\end{flushright}

\vspace{-2.truecm}
\vspace{1pt}
\begin{center}
{\Large \textbf{On the Convexity of Image of a\\[3pt] Multidimensional Quadratic Map}}
\end{center}
\vspace{6pt}
\begin{center}
{\large\textsl{Anatoly Dymarsky \\[10pt]}}
\textit{ Skolkovo Institute of Science and Technology,\\
Skolkovo Innovation Center, Nobel st.~3, Moscow, Russia, 143026\\ Department of Physics and Astronomy, University of Kentucky,\\ Lexington, KY 40506}\\ \vspace{6pt}
\end{center}

\vspace{6pt}
\begin{center}
\textbf{Abstract}
\end{center}
\vspace{-.1cm}
We study convexity of image of a general multidimensional quadratic map.  We split the full image into two parts by an appropriate hyperplane  such that one part is compact,  and formulate a sufficient condition for convexity of the compact part. We propose a way to identify such convex parts of the full image which can be used in practical applications. By shifting the hyperplane to infinity we extend the sufficient condition for convexity to apply to the full image of the quadratic map.  As a related result, we formulate a novel condition for the joint numerical range of $m$-tuple of hermitian matrices  to be convex. Finally, we illustrate our findings by considering several examples.  In particular we prove convexity of solvability set for the Power Flow equations in case of DC networks.

\vspace{6pt}
\begin{flushleft}
\end{flushleft}

\renewcommand{\thefootnote}{\arabic{footnote}}

\newpage
\tableofcontents

\section{Introduction and Main Results}
\label{intro}
In this paper we consider a multidimensional quadratic map of general form $f:\R^n\rightarrow \R^m$ (or $f:\C^n\rightarrow \R^m$) defined by  an $m$-tuple of symmetric (hermitian) matrices $A_i$, an $m$-tuple of vectors $v_i\in \R^n$ (or $v_i\in \C^n$), and a vector $f^0\in \R^m$,\footnote{\,Symbol ${}^*$ stands for transposition or hermitian conjugation depending on the context.}
\bea
\label{QM}
f_i(x)=x^* A_{i\,} x -v^*_i x-x^* v_i^{}+f^0_i\ .
\eea

An important question arising in many applications is when the image of the quadratic map \eqref{QM}, 
\bea
\label{fullimage}
{\mathscr F}(f)=\{f(x): x\in {\mathbb V}\} \subseteq \R^{m}\ ,\quad {\mathbb V}=\R^n {\rm\ or\ } {\mathbb V}=\C^n\ ,
\eea
is convex. This question was previously studied  in case of small $m=2,3$ or when only a few of $A_i$'s are linearly independent, see \cite{m2,specialquadraticconvex,Xia} and \cite{Polik} for references and a brief historic overview. Identifying necessary and sufficient conditions  for the convexity of  \eqref{fullimage} for general $m$ remains an open problem which we investigate in this paper.

Although the results concerning convexity of ${\mathscr F}(f)$ as a whole are scarce, something can be said about convexity of the image of $f$ locally. In particular, for any non-linear map $f$ there is an upper bound  on the size of a ball $B_{\varepsilon}(x_0)=\{x: |x-x_0|^2\le \varepsilon^2\}\subset {\mathbb V}$ such that its image $f(B_{\varepsilon}(x_0))=\{f(x): x\in B_{\varepsilon}(x_0)\}$ is convex \cite{Polyak}. Because of a very broad scope of this result the corresponding bound on $\varepsilon$ is always finite and $f(B_{\varepsilon}(x_0))$ can never fully cover ${\mathscr F}(f)$. The bound of \cite{Polyak} was improved in \cite{my} for the case of a general quadratic map \eqref{QM} and an  ellipsoid-shaped ball $B^+_{\varepsilon}(x_0)=\{x: |x-x_0|_+^2\le \varepsilon^2\}$ defined through some positive-definite matrix $A_+$, $|x|_+^2:=x^* A_+ x$,
\bea
\label{MY}
\varepsilon^{2}_{\rm max}(x_0,A_+):=\lim\limits_{\epsilon\rightarrow 0^+} \min\limits_{c\,\in\, {\mathcal C}} \left|(  c\cdot A-\lambda^+_{\rm min}(c\cdot A)A_++\epsilon)^{-1}c\cdot (v-A_{\,} x_0)\right|_+^2 \  ,\\
{\mathcal C}=\{c: c\in \R^m,\ |c|^2=1,\ \lambda^+_{\rm min}(c\cdot A)\le 0\}\ .
\eea
Above we introduced $\lambda^+_{\rm min}(A)$ to denote the smallest generalized eigenvalue of $A$ with respect to $A_+$.\footnote{The smallest generalized eigenvalue can be defines as $\lambda^+=\min_x {x^*A\, x\over  x^* A_{+}x}$.}  In \eqref{MY} and in what follows a sum of a matrix and a number understood in a sense that the number is multiplied by the identity matrix of an appropriate size. The dot product $\cdot$ stands for the standard scalar product in the Euclidean space.  Also notice the appearance of the ``plus" norm  in \eqref{MY}.

Reference \cite{my} proves that $f(B^+_{\varepsilon}(x_0))$ is strictly convex  for $\varepsilon\le \varepsilon_{\rm max}$. Now, if for  any non sign-definite  combination $c\cdot A \nsucc 0$, $c\in \R^m\backslash \{0\}$, the projection of the vector $c\cdot (v-A_{}x_0)$ on the eigenspace corresponding to the smallest  generalized  eigenvalue of $c\cdot A$ is non-vanishing, $\varepsilon_{\rm max}$ is infinite and the image $f(B^+_{\varepsilon}(x_0))$  is strictly convex for any $\varepsilon$. Consequently we have the following proposition. \\
{\bf Proposition 1.} If for any point $x_0\in {\mathbb V}$ and any $A_+\succ 0$ the value of $\varepsilon_{\rm max}(x_0,A_+)$ is infinite
 the full image ${\mathscr F}(f)$ is convex. This is a novel sufficient condition for the convexity of image of a general multidimensional quadratic map. 

The proof is trivial. For any two points $y_1,y_2\in {\mathscr F}(f)$ we  consider their pre-images $x_1, x_2\in {\mathbb V}$ (any pre-images if there are many). For a significantly large  $\varepsilon$, $x_1, x_2\in B^+_{\varepsilon}(x_0)$ and consequently all points $y(t)=y_1 (1-t)+y_2 t$ for $1\ge t \ge 0$ lie within $f(B^+_{\varepsilon}(x_0))\subset {\mathscr F}(f)$.

What happens if $\varepsilon_{\rm max}$ is very large but not infinite? Strictly speaking there is not much we can say about convexity of ${\mathscr F}(f)$ in this case. It is reasonable to expect that ${\mathscr F}(f)$ might still be very close to be convex because the image of a very large ball $B^+_{\varepsilon}(x_0)\subset{\mathbb V}$ is convex.  Unfortunately some points $x\in {\mathbb V},\ x\notin B^+_{\varepsilon}(x_0)$ might be mapped in $\R^m$ finite distance away from the origin  for any, even infinitely large $\varepsilon$ and this could spoil convexity of the full image.  If somehow we could arrange for $\thickmuskip=2mu \parallel y=f(x)\parallel$, $x\in {\mathbb V},\ x\notin B^+_{\varepsilon}(x_0)$, to be large with some appropriate norm $\parallel\ \parallel$ when $\varepsilon$ is large than not only ${\mathscr F}(f)$ would be convex when $\varepsilon^2_{\rm max}$ is infinite, but also we could outline a compact part of ${\mathscr F}(f)$ (an intersection of ${\mathscr F}(f)$ with a ``ball" of a certain size  defined by the norm $\parallel\ \parallel$),  which is convex  when  $\varepsilon^2_{\rm max}$ is finite. This idea, which we  develop in  section \ref{convexityFf}, leads to the following result.

For any vector $c_+\in \R^m\backslash \{0\}$ such that the combination $A_+:= c_+\cdot A$ is positive-definite, $A_+\succ 0$, let us define the   point $x_0=A_+^{-1}(c_+\cdot v)$. This is the unique point where the supporting hyperplane orthogonal to $c_+$ ``touches" ${\mathscr F}(f)$.
Let us introduce the following limit 
\bea
\label{MR4}
z_{\rm max}:=\lim\limits_{\epsilon\rightarrow 0^+} \min\limits_{c\,\in\, {\mathcal C}} \left|(  c\cdot A-\lambda_{\rm min}^+(c\cdot A)A_++\epsilon)^{-1}c\cdot (v- A_{\,} x_0)\right|_+^2 \  ,\\
{\mathcal C}=\{c: c\in \R^m,\ |c|^2=1,\ c\cdot c_+=0\}\ .
\eea

\noindent{\bf Proposition 2.} The compact set ${\mathscr F}(f,c_+,z_{\rm max})=\{f(x): x\in {\mathbb V},\ |x-x_0|^2_+\le  z_{\rm max}\}\subset {\mathscr F}(f)$ is convex. If $z_{\rm max}$ is infinite the whole image ${\mathscr F}(f)$ is convex. 

\noindent The Proposition 2 is reformulated in section \ref{zmax} in a more application-friendly way. See Proposition 2' and Comment 4 there. 

\begin{figure}
\vspace{3.5cm}
\begin{tikzpicture}
\node[] at (3,4) {};  
        \begin{subfigure}[b]{0.5\textwidth}
        \vspace{1.cm} 
                \includegraphics[width=1.\textwidth]{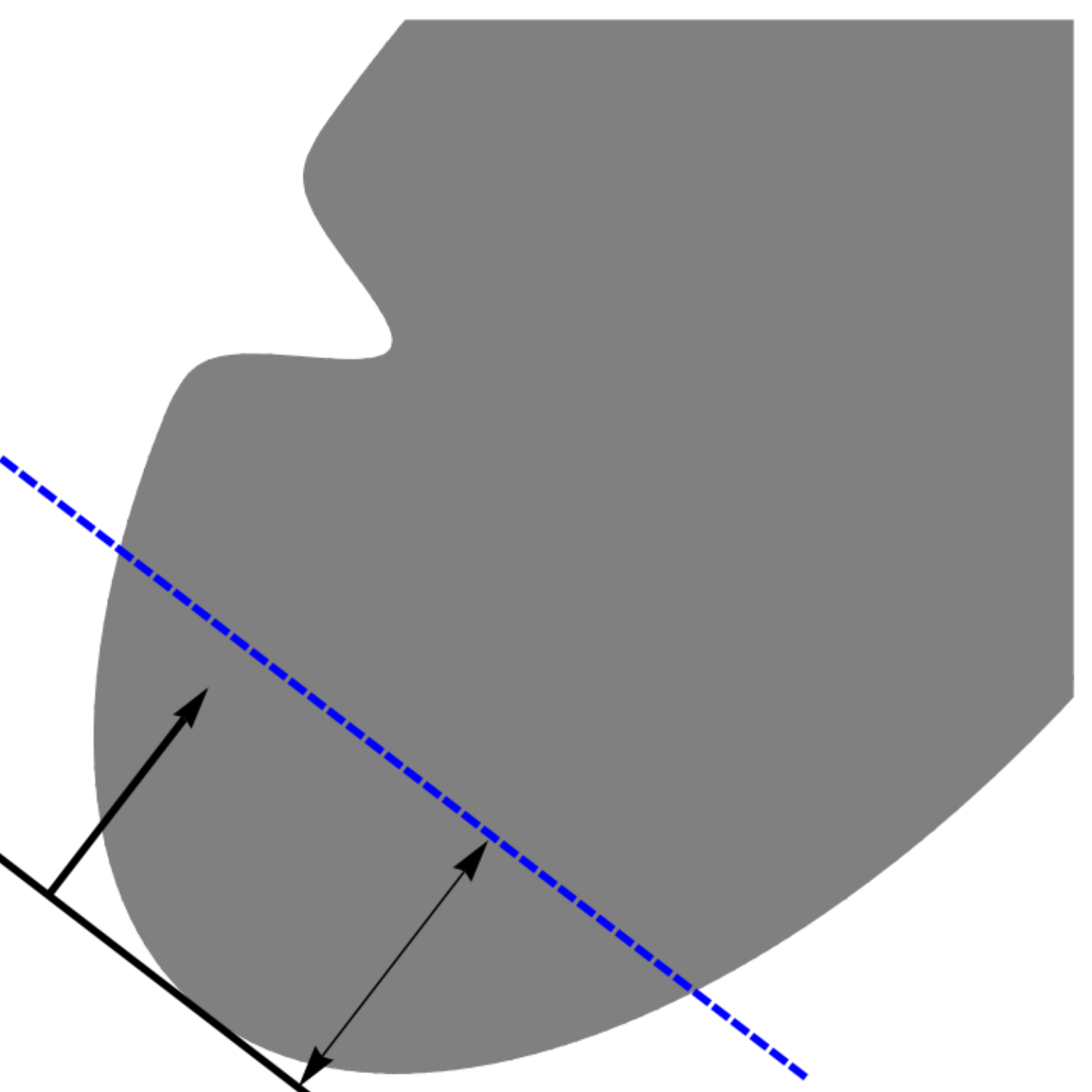}
                \caption{$a).$}
                \label{fig:idea}
                \vspace{-1.cm}
                \node at (.45,2.4) {$c_+$};
                \vspace{-1.cm}
                \node at (2.9,.65) {$z$};              
        \end{subfigure}~
        %
        \begin{subfigure}[b]{0.5\textwidth}
        \vspace{1.cm}
                \includegraphics[width=1.\textwidth]{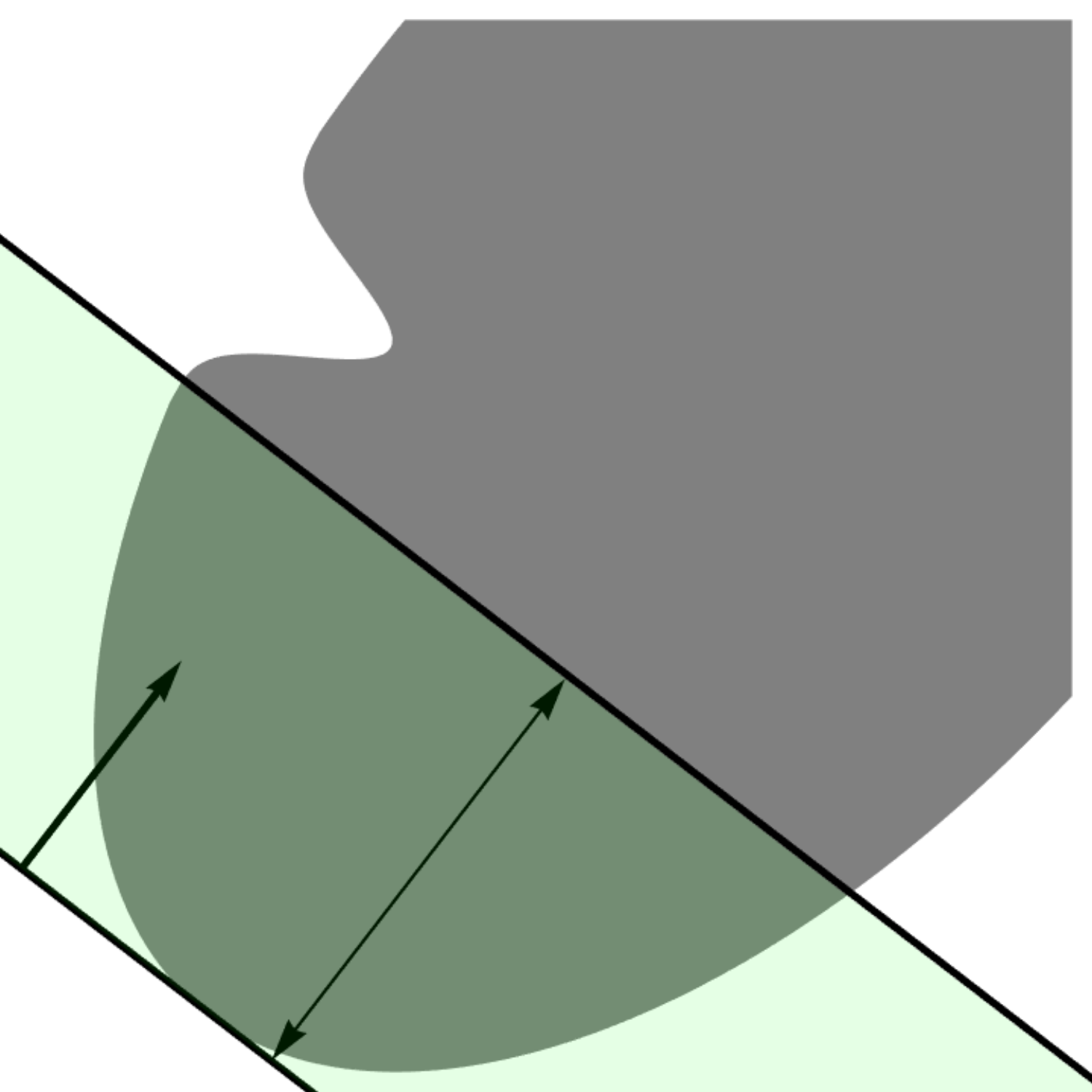}
                \caption{$\quad b).$}
                \label{fig:answer}
                \vspace{-1.cm}
                \node at (.3,2.5) {$c_+$};
                \vspace{-1.cm}
                \node at (3.4,1.15) {$z_{\rm max}$}; 
        \end{subfigure}        
\end{tikzpicture}
\vspace{1.cm}
\caption{(a) An illustration of the main idea. A hyperplane orthgonal to $c_+$ touches ${\mathscr F}(f)$ only at one point. If the boundary is smooth and strongly convex around that point a parallel hyperplane located a short distance apart will carve a compact subregion of ${\mathscr F}(f)$.  (b) $z_{\rm max}$ defines a maximal subregion of ${\mathscr F}(f)$ which is stably convex.}
\label{fig:plots}
\end{figure}
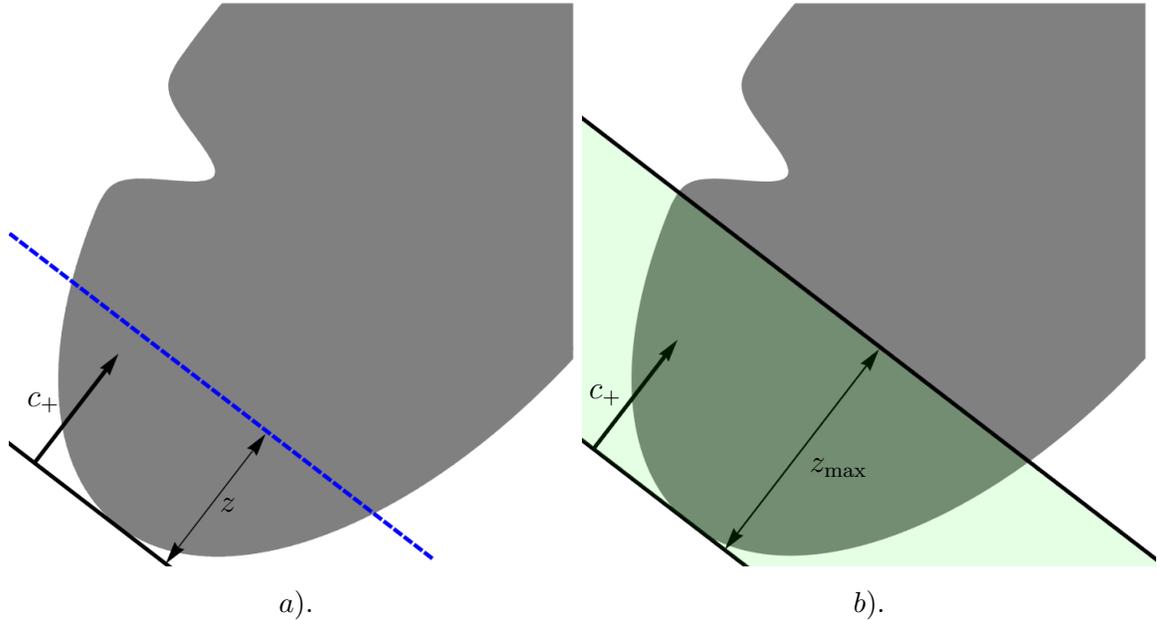

\noindent{\it Comment 1.} The convex set ${\mathscr F}(f,c_+,z_{\rm max})$ can be defined as a compact part of ${\mathscr F}(f)$ lying in a half-space defined by a certain hyperplane ${\mathscr H}_{c_+}(c_+\cdot f(x_0)+z_{\rm max})$, where for any $c\in \R^m\backslash\{0\}$ the corresponding orthogonal hyperplane is defined as 
\bea
 {\mathscr H}_{c}(F)=\{y: y\in \R^{m},\ c\cdot y=F\}\subset\R^{m}\ .
\eea
The hyperplane ${\mathscr H}_{c_+}(c_+\cdot f(x_0))$ is the supporting hyperplane of ${\mathscr F}(f)$ perpendicular to $c_+$. Hence ${\mathscr F}(f,c_+,z_{\rm max})$ is the part of ${\mathscr F}(f)$ bounded by two parallel hyperplanes separated by the distance $z_{\rm max}$ apart. This is illustrated in Fig.~\ref{fig:answer}.\\
{\it Comment 2.} The two criteria for convexity (Proposition 1 and Proposition 2) are complimentary to each other in the following sense. In \eqref{MY} $x_0$ is a regular point of $f$ and thus $\varepsilon_{\rm max}(x_0)$ may vanish. On the contrary in \eqref{MR4} $f(x_0)$ belongs to the boundary $\partial {\mathscr F}(f)$.\\
{\it Comment 3.} The sufficient condition of Proposition 1 depends on the choice of $x_0, A_+$. Potentially, if $\varepsilon_{\rm max}(x_0,A_+)$ is infinite for some  $x_0, A_+$ it might be finite for some other $x_0, A_+$. On the contrary if $z_{\rm max}$ \eqref{MR4} is infinite for some $c_+$ (and corresponding $x_0,A_+$), it will be infinite for all other choices of $c_+$ as well. If $z_{\rm max}$ is finite its value depends on $c_+$.
This is  explained in section~\ref{geomeaning}.

We derive the expression for $z_{\rm max}$ \eqref{MR4} and prove Proposition 2 in section~\ref{convexityFf}, while Proposition 1 directly follows from the results of \cite{my}. In section~\ref{connection} we discuss how our results pertaining to convexity of ${\mathscr F}(f)$ can be related to the classical question of convexity  of the joint numerical range of  $m$-tuple of symmetric (hermitian) matrices. In particular we formulate novel sufficient conditions for convexity of the joint numerical range in the subsection~\ref{newc}.
In section~\ref{zmax} we reformulate Proposition 2 and discuss different ways of calculating $z_{\rm max}$. In particular we outline an easy-to-calculate conservative  bound on $z_{\rm max}$  in the subsection~\ref{aprxsection}. In this way we formulate a practical way to outline a convex part of ${\mathscr F}(f)$, which can be used in applications.  Section 5 is devoted to applications and examples.  First, we prove convexity of solvability set for the Power Flow equations for the DC networks. Then we illustrate main points of this paper in case of several concrete  examples. In one case we consider an example when the image is not convex and calculate $z_{\rm max}$ both exactly and approximately, using the bound of the subsection~\ref{aprxsection}.

\section{Convexity of ${\mathscr F}(f)$} 
\label{convexityFf}
Let us first approach the question of convexity of ${\mathscr F}(f)$ locally. From now on we assume that the $m$-tuple of matrices $A_i$ is {\it definite} in the sense of reference \cite{Sheriff}, i.e.~there is a positive-definite combination $A_+:=c_{+}\cdot A\succ 0$ for some $c_+\in \R^m\backslash\{0\}$. The corresponding supporting hyperplane ${\mathscr H}_{c_+}(z_0)$ touches ${\mathscr F}(f)$ only at one point $y_0=f(x_0)$,
\bea
\label{cplus}
 x_0=A_+^{-1\,}v_+\ ,\quad z_0=-v_+^*A_+^{-1\,}v_+\ , \quad v_+:=c_+\cdot v\ .
\eea
This is schematically depicted in Fig.~\ref{fig:idea}. Provided the boundary of ${\mathscr F}(f)$ is smooth and strongly convex at $y_0$ it is tempting to say that the compact part of ${\mathscr F}(f)$ lying in a half-space defined by the hyperplane ${\mathscr H}_{c_+}(z_0+z),\ z> 0$, would be convex, at least for very small $z$. 
It is easy to see that the compact part of ${\mathscr F}(f)$ bounded by ${\mathscr H}_{c_+}(z_0+z),\ z> 0$,
is the image of the ball $B^+_{z}(x_0)=\{x: |x-x_0|^2_+\le z\}$ under the map ${f}$. We would like to find an upper bound on $z$ such that  the image ${f}(B^+_{z}(x_0))$ is convex.  This problem almost identically repeats the question studied by Polyak \cite{Polyak} for general ${f}$, and further investigated in \cite{my} in case of quadratic $f$, with one crucial distinction: point $x_0$ is not a regular point of $ f$. 

The following observation drastically simplifies further analysis: each ellipsoid $|x-x_0|_+^2=z$ is mapped into its own hyperplane ${\mathscr H}_{c_+}(z_0+z)$. Hence the convexity of ${f}(B^+_{z}(x_0))$ requires convexity of ${f}(|x|_+^2=z')$ for all $z'\le z$. Up to a translation and a trivial change of basis in ${\mathbb V}$ the image ${f}(|x-x_0|^2_+=z)$ is nothing but the inhomogeneous joint numerical range, the notion we introduced in \cite{my},
\bea
\label{gjnr}
{\mathcal F}(\bf A,v)&=&\left\{y_i: \exists\ x,\ y_i=x^* {\bf A}_{i\,} x -{\bf v}^*_i x-x^* {\bf v}_i^{},\ |x|^2=1\right\}\ ,
\eea
for an $m$-tuple of symmetric (hermitian) matrices $\bf A_i$ and a $m$-tuple of vectors $\bf v_i$.
There it was proven that if
\bea
\label{suff}
\lim\limits_{\epsilon\rightarrow 0^+}\min_{c\neq 0} \left| \left(c\cdot {\bf A} -\lambda_{\rm min}(c\cdot {\bf A})+\epsilon\right)^{-1}c\cdot {\bf v}\right| \ge 1\ ,
\eea
${\mathcal F}({\bf A,v})$ is strictly convex (strongly convex for strong inequality in \eqref{suff}) and smooth. In fact \eqref{suff} is a criterion for stable convexity, i.e.~impossibility to ruin convexity of ${\mathcal F}({\bf A,v})$ by an infinitesimal deformation of $\bf A,v$.

Applying this directly to the image $f(|x-x_0|^2_+=z)$ will fail because the latter  is ``flat'', i.e.~it lies within the hyperplane ${\mathscr H}_{c_+}(z_0+z)$ and hence can not be strictly convex. This is easy to fix by considering an orthogonal  projection $\varphi$ of  $\R^m$ on ${\mathscr H}_{c_+}\simeq \R^{m-1}$. Now the criterion \eqref{suff} can be applied directly to the $(m-1)$-dimensional inhomogeneous joint numerical range $(\varphi\circ f)(|x-x_0|^2_+=z)$, 
\bea
\label{suff2}
z_{\rm max}:=\lim\limits_{\epsilon\rightarrow 0^+}\min_{{ c}\neq 0} \left| \left({ c}\cdot {A} -\lambda^+_{\rm min}({ c}\cdot {A})A_+ +\epsilon\right)^{-1}{ c}\cdot (v-A_{\, }x_0)\right|_+^2 \ge z\ .
\eea
The minimum \eqref{suff2} is taken over the space of the equivalence classes $c \in \R^{m},\ c\simeq c+\mu_{} c_+,\ \forall \mu \in \R$. Since the minimized expression is homogeneous in $c$ the condition $c\neq 0$ could be substituted by $|c|^2=1,\ c\cdot c_+=0$. (Notice, that if $c\propto c_+$, \eqref{suff2} vanishes.)

Clearly, convexity of ${\mathscr F}(f)\cap {\mathscr H}_{c_+}(z_0+z)$ for all $z\le z_{\rm max}$ is a necessary condition for the convexity of $f(B^+_{z_{\rm max}}(x_0))$ but may  not be sufficient. To establish convexity of $f(B^+_{z_{\rm max}}(x_0))$ let us choose an arbitrary vector $c\in \R^m\backslash\{0\}$ and find an intersection of $f(B^+_{z_{\rm max}}(x_0))$ with a supporting hyperplane ${\mathscr H}_{c}(F_c)$ orthogonal to $c$ which touches $f(B^+_{z_{\rm max}}(x_0))$ ``from below'',
\bea
F_c=\min_{x\in B^+_{z_{\rm max}}(x_0)} c\cdot f(x)\ .
\eea
For $c$ collinear with $c_+$ the answer is simple. When $c=\alpha \, c_+$, $c\cdot c_+>0$, the hyperplane ${\mathscr H}_{c_+}(z_0)$ intersects $f(B^+_{z_{\rm max}}(x_0))$ at the unique point $y_0$. When $c\cdot c_+<0$ the intersection ${\mathscr H}_{c_+}(z_0+z_{\rm max})\cap f(B^+_{z_{\rm max}}(x_0))$ is a $(m-1)$-dimensional convex set $f(|x-x_0|_+^2=z_{\rm max})$. 
For $c$ not collinear with $c_+$ we can first find a conditional minimum for $|x-x_0|_+^2=z$ and then minimize with respect to $z$ in the interval $[0,z_{\rm max}]$. This problem was solved in \cite{my} where it was shown that the minimum is achieved at such $z$ that $\lambda(z)$, uniquely determined by the conditions
\bea
\left| \left(c\cdot {A} -\lambda(z) A_+ \right)^{-1}c\cdot (v-A_{\, }x_0)\right|_+^2 = z\ ,\quad \lambda(z)\le \lambda^+_{\rm min}(c\cdot A)\ ,
\eea
is equal to $\min\{0,\lambda(z_{\rm max})\}$. Since $\lambda(z)$ is a monotonically increasing function of $z$ on the interval $[0,z_{\rm max}]$ such a point is unique. Correspondingly the supporting hyperplane ${\mathscr H}_{c}(F_c)$ intersects $f(B^+_{z_{\rm max}}(x_0))$ at a unique point.  Hence $\partial_{\,}{\rm Conv}[f(B^+_{z_{\rm max}}(x_0))]\subset f(B^+_{z_{\rm max}}(x_0))$ where $\rm Conv$ stands for convex hull. Now, to prove that $f(B^+_{z_{\rm max}}(x_0))$ is convex we need to show that  any point of ${\rm Conv}[f(B^+_{z_{\rm max}}(x_0))]$  belongs to $f(B^+_{z_{\rm max}}(x_0))$. This is obviously true because an intersection of $f(B^+_{z_{\rm max}}(x_0))$ with the hyperplane ${\mathscr H}_{c_+}(F)$ for any $F$ is either empty or convex. This finishes the proof of Proposition 2. 

\subsection{Geometrical meaning of $z_{\rm max}$}
\label{geomeaning}
To interpret $z_{\rm max}$ geometrically we would need to understand different scenarios of how ${\mathscr F}(f)$ may intersect with its supporting hyperplanes. Let us consider a vector $c\in \R^m\backslash\{0\}$ and find a supporting hyperplane to ${\mathscr F}(f)$ that is orthogonal to $c$. There are several possible scenarios. First, $c\cdot A$ is sign-definite. The corresponding supporting hyperplane intersects  ${\mathscr F}(f)$ at the unique point $f(x),\ x=(c\cdot A)^{-1}c\cdot v$. Second, $c\cdot A$ has both positive and negative eigenvalues. There is no corresponding supporting hyperplane in this case because  ${\mathscr F}(f)$ stretches to infinity in both directions along $c$. Finally, $c\cdot A$ is semi-definite and degenerate. There are two possibilities (Fredholm alternative) in this case. If the equation $(c\cdot A)_{}x=c\cdot v$ admits no solution, there is no supporting hyperplane to ${\mathscr F}(f)$ orthogonal to $c$ because ${\mathscr F}(f)$ stretches to infinity in both directions along $c$. Another option is when there is a whole linear space of solutions of $(c\cdot A)x=c\cdot v$. Each such solution $x$ corresponds to a point $f(x)$ from the boundary $\partial{\mathscr F}(f)$ belonging to the same  supporting hyperplane orthogonal to $c$. 

In case a supporting hyperplane intersects ${\mathscr F}(f)$ over more than one point we would like to call all points of $\partial{\mathscr F}(f)$ belonging to this supporting hyperplane a ``flat edge''.\footnote{More precisely we should define ``flat edge'' as an intersection of ${\mathscr F}(f)$ with a supporting hyperplane if the pre-image of this intersection consists of more than one point.  But a ``flat edge'' defined this way can consist of one point only if all matrices $A_i$ have a common zero eigenvector which is in contradiction with the assumption that the set of $A_i$'s is definite.} If $\partial {\mathscr F}(f)$ includes ``flat edges'', ${\mathscr F}(f)$ can not be strictly convex  and certainly ${\mathscr F}(f)$ can not be stably convex. In general such ${\mathscr F}(f)$ will not be convex at all.

Now let us look at the definition of $z_{\rm max}$ \eqref{MR4}.
For any $c$, $c\cdot c_+=0$, the limit $\epsilon\rightarrow 0^+$ in \eqref{MR4} will be finite only if the equation  
\bea
\label{cdef}
({\bf c}\cdot A) (x+x_0)={\bf c}\cdot v\ ,\quad {\bf c}(c)=c-\lambda^+_{\rm min}(c\cdot A)c_+\ ,
\eea has nontrivial solution(s). Since ${\bf c}\cdot A$ is degenerate and positive semi-definite, existence of nontrivial solutions is the same as the existence of a ``flat edge'' orthogonal to ${\bf c}$. Hence $z_{\rm max}$ will be infinite unless there is a supporting hyperplane touching ${\mathscr F}(f)$ at more than one point. The latter is the property of $f$ and ${\mathscr F}(f)$ and does not depend on the choice of $c_+$ and $x_0,A_+$. 

\begin{figure}
\vspace{1.cm}
\begin{tikzpicture}
        \begin{subfigure}[b]{0.5\textwidth}
        \vspace{1.cm} 
                \includegraphics[width=1.\textwidth]{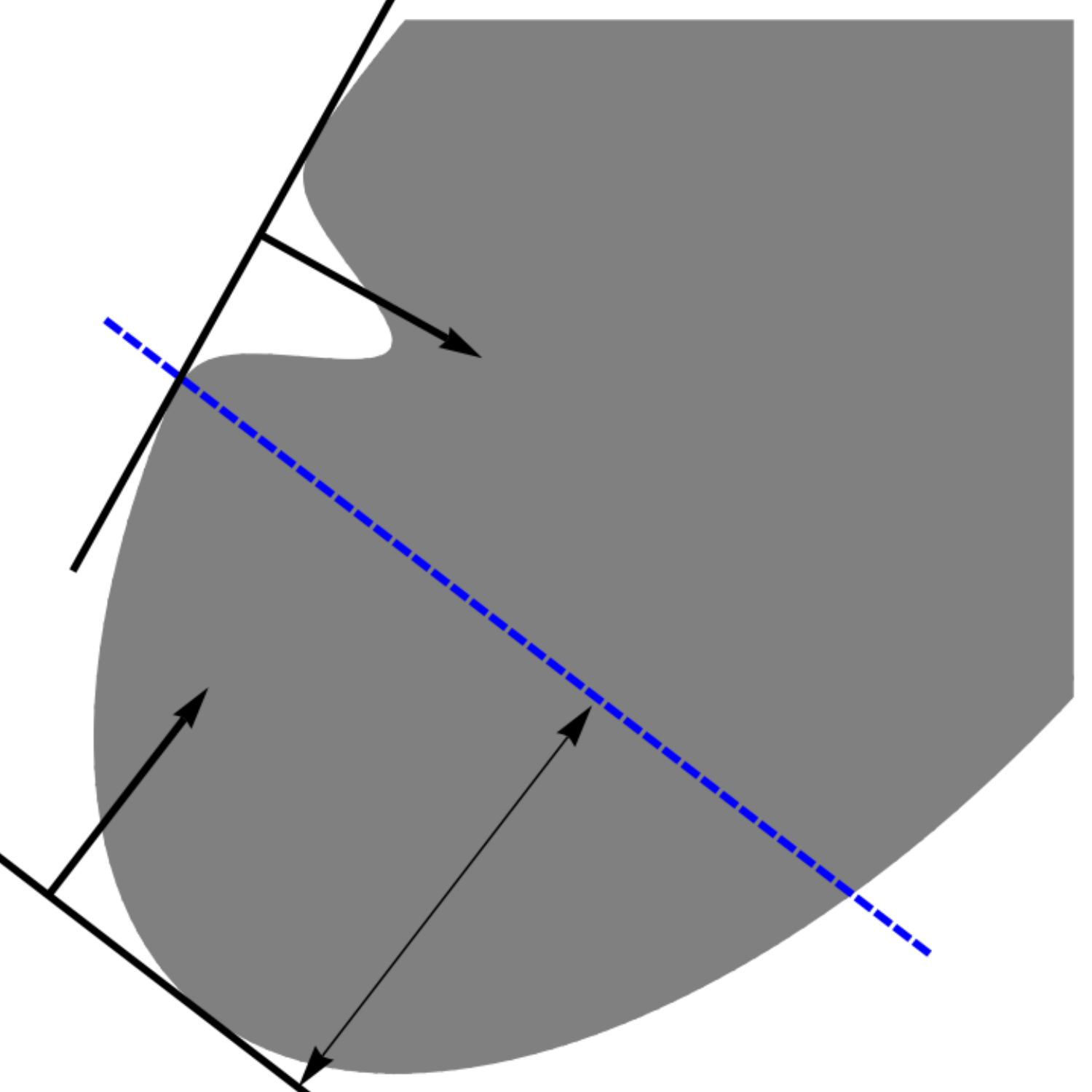}
                \caption{$a).$}
                \label{fig:meaningz}
                \vspace{-1.cm}
                \node at (.46,2.5) {$c_+$};
                \vspace{-1.cm}
                \node at (3.65,1.15) {$z_{\rm max}$};    
                \node at (3.4,6) {$c$};             
        \end{subfigure}~ 
        %
        \begin{subfigure}[b]{0.5\textwidth}
            \label{fig:noncon}      
                \includegraphics[width=1.\textwidth]{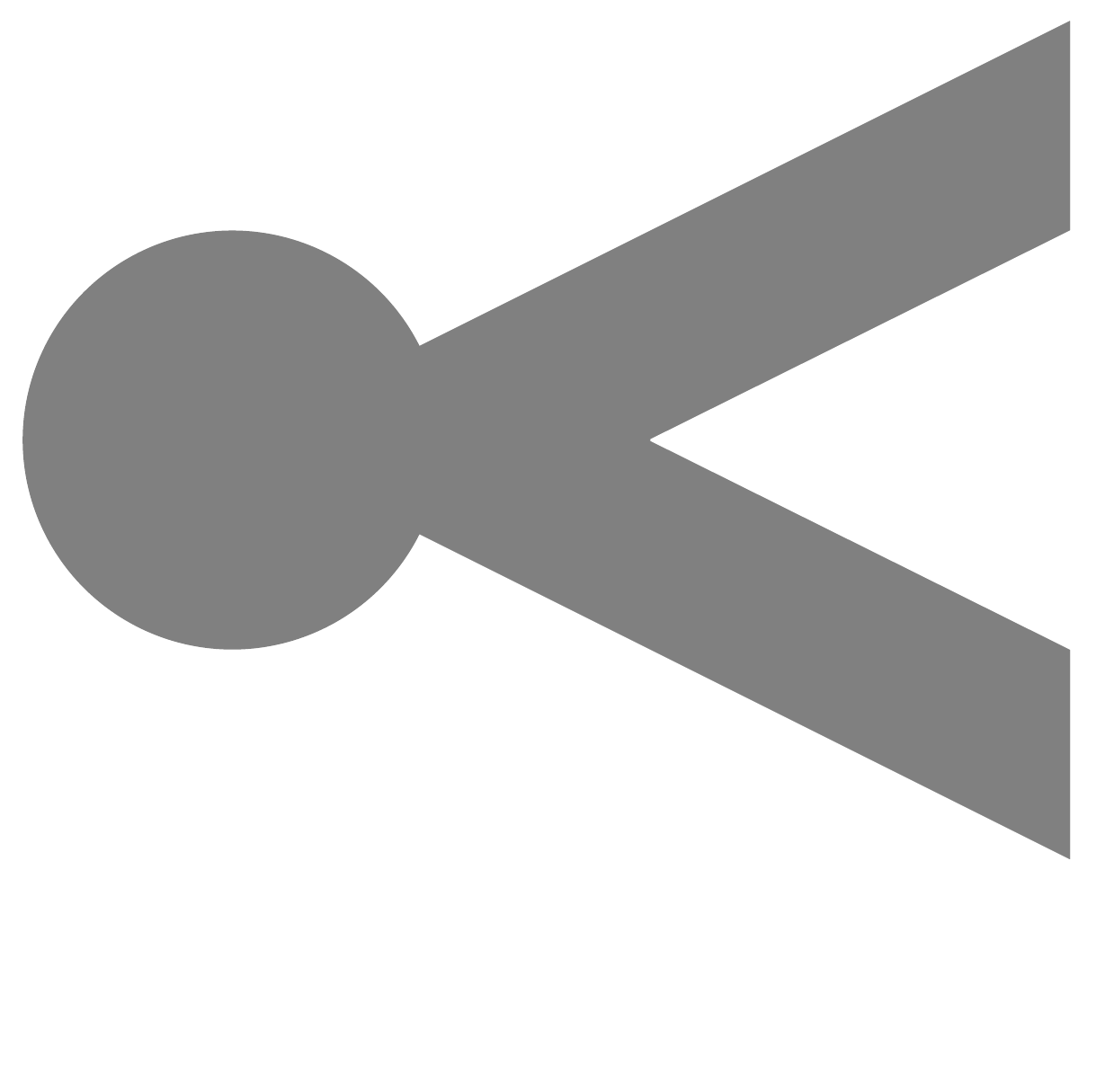}  \\
                \vspace{-.08truecm}            
                 \node at (4.,0){$b).$};
        \end{subfigure}        
\end{tikzpicture}

\caption{(a) Geometrical meaning of $z_{\rm max}$. It is the distance from a supporting hyperplane orthogonal to $c_+$ to the closest point  belonging to a``flat edge''.  (b) A non-convex figure which intersects any supporting hyperplane at exactly one point. (The shevron-shape strips continue to infinity.)}
\label{fig:plots}
\end{figure}

It is easy to show that for any  vector $c$ such that ${\bf c}(c)$ is orthogonal to a ``flat edge'' (namely, $({\bf c}\cdot A)\succeq 0,$ but $({\bf c}\cdot A)\nsucc 0$ and the space of solutions $ {\mathcal S}=\{x: ({\bf c}\cdot A)x={\bf c}\cdot v\}$ is non-trivial), the limit $\epsilon\rightarrow 0^+$ in \eqref{MR4} calculates  the minimum value of $c_+\cdot f(x),\ x\in {\mathcal S}$. 
Consequently $z_{\max}$ is the distance from the supporting hyperplane orthogonal to  $c_+$ to a closest point $y\in\partial{\mathscr F}(f)$ belonging to any ``flat edge'' inside $\partial{\mathscr F}(f)$. This is depicted in Fig.~\ref{fig:meaningz}. This also explains that ${\mathscr F}(f,c_+,z)$ is stably convex for $z\le z_{\rm max}$ and is not stably convex $z>z_{\rm max}$. 

A comment is in order. If ${\mathscr F}$ were compact, absence of ``flat edges''  would immediately guarantee that the outer boundary $\partial {\mathscr F}$ is a boundary of the convex hull of ${\mathscr F}$ while the latter is strictly convex. But in a non-compact case this is not true. Fig.~2b provides a non-compact example without ``flat edges'' (meaning that each supporting hyperplane is touching the figure is exactly at one point) when the outer boundary is not the boundary of the convex hull. Hence the proof that the absence of ``flat edges'' ($z_{\rm max}\rightarrow \infty$) implies that the  boundary $\partial {\mathscr F}$ confines a convex set, which was presented in section~\ref{convexityFf},  was not superfluous. 

\section{Connection with the joint numerical range}
\label{connection}
In this section we would like to look at the convexity of ${\mathscr F}(f)$ from a slightly different angle. Shifting $f^0_i$ by a constant does not affect convexity and therefore without loss of generality we assume $f^0_i=0$. Then ${\mathscr F}(f)$ can be thought of as an intersection of the image of the ``extended'' {\it homogeneous} quadratic map ${\mathsf f}:\R^{n+1}\rightarrow \R^{m+1}$ (or ${\mathsf f}:\C^{n+1}\rightarrow \R^{m+1}$) and a hyperplane
\bea
\label{FFH}
&&{\mathscr F}(f)={\mathscr F}({\mathsf f}) \cap {\mathscr H}_{e_{m+1}}(1)\ , 
\eea
where $e_{m+1}$ is the $(m+1)$-th  basis vector and
\bea
\label{qm}
&&{\mathsf f}_I={\mathsf x}^* {\mathsf A}_I {\mathsf x}\ , \ \quad \qquad I=1,\dots {m+1}\ ,\qquad\quad \ \ \,  {\mathsf x}\in \R^{n+1}\ ({\rm or\ } \C^{n+1})\ ,\\ 
\label{matA} &&{\mathsf A}_i=\left(
\begin{array}{c|c}
A_i & -v_i \\ \hline 
-v_i^* & 0
\end{array}
\right),\ \ i=1,\dots,m\ , \quad {\mathsf A}_{m+1}=\left(
\begin{array}{c|c}
0^{}_{n\times n} & 0_{n\times 1} \\ \hline
0_{1\times n} & 1_{}^{}
\end{array}
\right)\ ,\\
&& {\mathscr H}_{{\mathsf c}}(F)=\{y: y\in \R^{m+1}\ ,\ {\mathsf c}\cdot y=F\}\subset\R^{m+1}\ ,\ \  \, \quad \quad \,  \forall\ {\mathsf c}\in \R^{m+1}\ .
\eea
Obviously, convexity of ${\mathscr F}(\mathsf f)$ would imply convexity of ${\mathscr F}(f)$. Validity of the converse statement is a subject of the following discussion. On general grounds convexity of ${\mathscr F}(f)$ together with the conic structure of ${\mathscr F}({\mathsf f})$ and the relation \eqref{FFH} do not imply convexity of ${\mathscr F}({\mathsf f})$  but only of ${\mathscr F}({\mathsf f}) \backslash {\mathscr H}_{e_{m+1}}(0)$.

Possible advantage of introducing ${\mathscr F}(\mathsf f)$ is that 
 convexity of homogeneous quadratic maps was studied previously in \cite{Sheriff}. There is a general relation between definite homogeneous quadratic maps and joint numerical ranges of some auxiliary matrices.
Geometrically ${\mathscr F}(\mathsf f)$ is  a cone. If $\mathsf f$ is definite, i.e.~there is a linear combination ${\mathsf c}_+\cdot {\mathsf A}\succ 0$, the base of the cone ${\mathscr F}(\mathsf f)\cap {\mathscr H}_{{\mathsf c}_+}(F),\ F>0$, is compact and equal (up to a linear isomorphism) to 
a joint numerical range ${\mathcal F}({\mathcal A})$ of $m$ matrices ${\mathcal A}_i$ which can be explicitly constructed from ${\mathsf A}_I$.
Hence proving convexity of ${\mathscr F}(\mathsf f)$ is the same as establishing convexity of  ${\mathcal F}({\mathcal A})$. The latter is a question with a long and rich history. 

For an $m$-tuple of symmetric or hermitian matrices ${ \mathcal A_i}$ the joint numerical range is defined as follows
\bea
\label{JNR}
{\mathcal F}(\mathcal A)=\left\{y_i: \exists\ { x},\  x\in {\mathbb V},\ |x|^2=1,\ y_i=x^* {\mathcal A}_{i\,} x \right\} \subset \R^m\ ,
\eea
where $\mathcal A_{i\,}$ are symmetric or Hermitian matrices and  ${\mathbb V}=\R^n$ or ${\mathbb V}=\C^n$ correspondingly. The question of convexity of ${\mathcal F}$ goes back to Housdorff and Toeplitz \cite{HT} who proved  that ${\mathcal F}$ is always convex for $m=2$ hermitian matrices (i.e.~${\mathbb V}=\C^n$) and $n>1$. There are numerous results for small $m=2,3$ \cite{Dines, Brickman,Barvinok,vectorv} (also see \cite{Poon} for references) and a few specific conditions rendering 
${\mathcal F}$ non-convex. The case of general $m,n$ is not fully understood, although there is a sufficient condition that guarantees that ${\mathcal F}(\mathcal A)$ is strongly convex and smooth: if  for all linear combinations of ${\mathcal A}_i$ the dimension of the eigenspace corresponding to the lowest eigenvalue is the same \cite{Gutkin}. In terms of the corresponding map ${\mathsf f}$ this is the condition that for any linear combination such that ${\mathsf c}\cdot  {\mathsf A}\succeq 0$, but ${\mathsf c}\cdot  {\mathsf A}\nsucc 0$,  dimension of ${\rm ker}({\mathsf c}\cdot  {\mathsf A})$ is the same. This condition is a generalization of {\it roundness} defined in \cite{Sheriff} which guarantees ``roundness'' (strict convexity and smoothness) of the base of the cone ${\mathscr F}({\mathsf f})$.

We will see later that in our case of interest \eqref{matA} the condition that ${\rm dim(ker}(c\cdot  {\mathsf A}))$ remains the same for different ${\mathsf c}$ is not satisfied. Thus the results of \cite{Gutkin}, \cite{Sheriff} do not help to establish convexity of ${\mathscr F}(f)$. In the following we will reverse the logic and extend the sufficient condition for convexity of ${\mathscr F}(f)$, Proposition 2, to  ${\mathscr F}({\mathsf f})$ and $\mathcal F({\mathcal A})$. In this way we formulate new criteria for the convexity of the joint numerical range.

From now on we assume that $z_{\rm max}$ given by \eqref{MR4} is infinite, which implicitly assumes  the set of $A_i$'s is definite. 
Let us show that  the set of ${\mathsf A}_I$'s given by \eqref{matA} is definite as well. Starting from an appropriate $c_+$, let us consider the vector ${\mathsf c}_+\in \R^{m+1}$,
\bea
\label{cc}
{\mathsf c}_{+i}=c_{+i},\ i=1,\dots,m,\qquad {\mathsf c}_{+(m+1)}>v_+^* A_+^{-1} v_+\ .
\eea
Then the matrix ${\mathsf A}_+:={\mathsf c}_+\cdot {\mathsf A}$ is positive-definite as follows from the Sylvester's criterion.

\begin{SCfigure}[50]
\captionsetup{singlelinecheck=off}
    \caption[.]{Join numerical range of two matrices (see section \ref{3bus} for details)
    \begin{displaymath}
    {\mathcal A}_1=\left(
\begin{array}{ccc}
 0 &  \sqrt{\frac{8}{15}} & \frac{1}{\sqrt{10}} \\
  \sqrt{\frac{8}{15}} & 0 & 0 \\
 \frac{1}{\sqrt{10}} & 0 & 0 \\
\end{array}
\right), \  {\mathcal A}_2=\left(
\begin{array}{ccc}
  0 & 0 & 0\\
    0 & 0 & 0\\
      0 & 0 & 1
  \end{array}
\right).
    \end{displaymath}}
\includegraphics[width=0.4\textwidth]{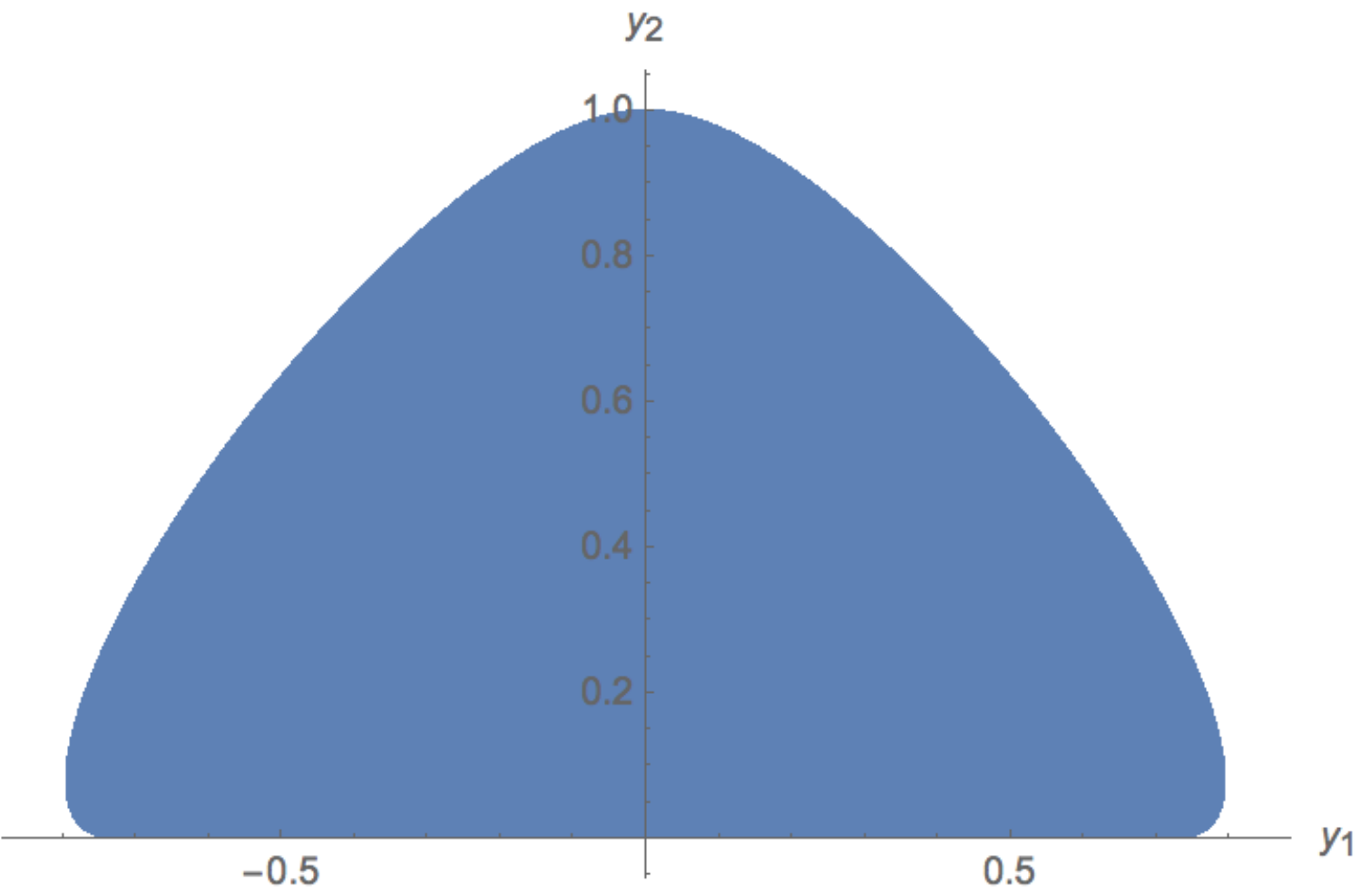}
\label{fig:flatedge}
\end{SCfigure}

Infinite value of  $z_{\rm max}$ implies  that for any $c$, $c\cdot A\succeq 0,\ c\cdot A\nsucc 0$, equation $(c\cdot A)x=c\cdot v$ has no solution (see section~\ref{geomeaning}). 
This immediately implies that any ${\mathsf x}\in {\ker}({\mathsf c}\cdot {\mathsf A})$, where ${\mathsf c}\cdot {\mathsf A}\succeq 0$ is either trivial or must satisfy ${\mathsf x}^*{\mathsf A}_{m+1} {\mathsf x}>0$, unless ${\mathsf c}\propto e_{m+1}$. Hence ${\dim}({\ker}({\mathsf c}\cdot {\mathsf A}))=1$ for all appropriate ${\mathsf c}$, except for the special case ${\mathsf c}\propto e_{m+1}$ when  ${\dim}({\ker}({\mathsf A}_{m+1}))=n$. Presence of this exceptional direction makes it impossible to apply the results of  \cite{Gutkin, Sheriff}. Indeed those works were focused on strictly convex and smooth joined numerical range ${\mathcal F}({\mathcal A})$, properties guaranteed by constancy of
${\dim}({\ker}({\mathsf c}\cdot {\mathsf A}))=1$ for all appropriate ${\mathsf c}$ such that  ${\mathsf c}\cdot {\mathsf A}\succeq 0$, but ${\mathsf c}\cdot {\mathsf A}\nsucc 0$. On the contrary, in our case 
${\mathcal F}({\mathcal A})$ has a ``flat edge'' perpendicular to the direction of $e_{m+1}$ (after projection on ${\mathscr H}_{{\mathsf c}_+}(1)$), which is a direct consequence of ${\dim}({\ker}({\mathsf c}\cdot {\mathsf A}))=n$ for one particular direction ${\mathsf c}\propto e_{m+1}$. The joint numerical range of this kind -- with one side being flat (flat edge) -- is shown in Fig.~\ref{fig:flatedge}.

The ``flat edge'' of ${\mathcal F}(\mathcal A)$ is an $(m-1)$-dimensional figure. Let us first establish it is convex.  We introduce 
\bea
{\mathfrak F}(F)={\mathscr F}({\mathsf f})\cap  {\mathscr H}_{{\mathsf c}_+}(1)\cap {\mathscr H}_{e_{m+1}}(F),\ F\ge 0 \,
\eea
and notice that  ${\mathfrak F}(0)$ is isomorphic to the ``flat edge'' of ${\mathcal F}(\mathcal A)$.
For any $F>0$, ${\mathfrak F}(F)$ is isomorphic to $f(|x-x_0|_+^2=z)={\mathscr F}(f)\cap {\mathscr H}_{c_+}(z_0+z)$ where corresponding $c_+$ (and hence $x_0,z_0$) is related to ${\mathsf c}_+$ through \eqref{cc} and $z$ is some function of $F$ an other parameters. The limit $F\rightarrow 0$ corresponds to $z\rightarrow \infty$. For any sufficiently small $F_0>F>0$, where $z(F_0)=z_0$,  we proved that ${\mathfrak F}(F)=f(|x-x_0|_+^2=z)$ is strongly convex. Hence by continuity ${\mathfrak F}(0)$ is convex as well.

The rest of the proof closely follows the logic outlined in section~\ref{convexityFf}. Let us consider a vector ${\mathsf c}\in\R^{m+1}$, $|{\mathsf c}|^2=1$, ${\mathsf c}\cdot {\mathsf c}_+=0$, and find an intersection of ${\mathcal F}(\mathcal A)\simeq {\mathscr F}({\mathsf f})\cap  {\mathscr H}_{{\mathsf c}_+}(1)$ with the supporting hyperplane ${\mathscr H}_{{\mathsf c}}(F_{\mathsf c})$ orthogonal to $\mathsf c$  touching  ${\mathcal F}(\mathcal A)$ ``from below''
\bea
F_{\mathsf c}=\min\limits_{y\in  {\mathscr F}({\mathsf f})\cap  {\mathscr H}_{{\mathsf c}_+}(1)} {\mathsf c}\cdot y\ .
\eea
For any ${\mathsf c}$ the intersection consists of a unique point except  for ${\mathsf c}=e_{m+1}$ when it is a convex ``flat edge''. Hence the boundary $\partial {\mathcal F}(\mathcal A)$ is the boundary of the convex hull of ${\mathcal F}(\mathcal A)$, $\partial_{\,}{\rm Conv}[ {\mathcal F}(\mathcal A)]\subset {\mathcal F}(\mathcal A)$. Finally, since the intersection ${\mathscr F}({\mathsf f})\cap  {\mathscr H}_{{\mathsf c}_+}(1)\cap {\mathscr H}_{e_{m+1}}(F)$ is convex (or empty) for any $F$ we conclude that all points confined by $\partial {\mathcal F}(\mathcal A)$ belong to ${\mathcal F}(\mathcal A)$. This establishes convexity of ${\mathscr F}({\mathsf f})$ and $\mathcal F({\mathcal A})$ provided the criterion for convexity of ${\mathscr F}(f)$, Proposition 2, is satisfied.

\subsection{New Criteria for Convexity for the Joint Numerical Range}
\label{newc}
The proof of convexity of ${\mathcal F}(\mathcal A)$ can be cast in a form of a self-contained criterion for convexity based only on the properties of matrices ${\mathcal A}_i$. 

Let us consider a joint numerical range \eqref{JNR} defined by an $m$-tuple of $n\times n$ symmetric (hermitian) matrices $\mathcal A_i$. If for any linear combination $c\cdot \mathcal A$, $c\in \R^m\backslash\{0\}$ its smallest eigenvalue is not degenerate, ${\mathcal F}(\mathcal A)$ is convex \cite{Gutkin}. If for any linear combination $c\cdot \mathcal A$, $c\in \R^m\backslash\{c: c=\mu\, e,\ \mu\ge 0\}$, where $e\in \R^m\backslash\{0\}$ is a fixed vector, its smallest eigenvalue is not degenerate, but  the smallest eigenvalue of $e\cdot \mathcal A$ is $(n-1)$-times degenerate, ${\mathcal F}(\mathcal A)$ is convex if the following condition is satisfied. We add the identity matrix ${\mathbb I}_{n\times n}$ to the set of $\mathcal A_i$'s thus bringing the total number of matrices to $(m+1)$. By changing a basis in ${\mathbb V}$ we bring $e\cdot \mathcal A$ to the form of ${\mathsf A}_{m+1}$ from \eqref{matA}. By taking a linear combination of $e\cdot \mathcal A$ with other $m$ matrices we bring them to the form of ${\mathsf A}_i$ from \eqref{matA} and in this way define $m$ vectors $v_i$ and $(n-1)\times (n-1)$ matrices $A_i$. Now, if the auxiliary map $f$ defined by $A_i, v_i$ through \eqref{QM} satisfies the convexity criterion of Proposition 2, $z_{\rm max}\rightarrow \infty$, then ${\mathcal F}(\mathcal A)$ is convex. This is the {\it first} new sufficient condition for convexity. 

In fact there is another sufficient criterion of convexity of ${\mathcal F}(\mathcal A)$ ``buried'' inside the proof in section~\ref{connection}. Indeed, the ``flat edge'' ${\mathfrak F}(0)$ is the joint numerical range associated with the homogeneous quadratic map $y_i=x^* A_i x$. Hence the {\it second} new sufficient condition for the convexity of the joint numerical range can be formulated as follows.
Let us consider a joint numerical range \eqref{JNR} defined by an $m$-tuple of $n\times n$ symmetric (hermitian) matrices $\mathcal A_i$. By adding the identity matrix to the set of $\mathcal A_i$'s we define an $(m+1)$-tuple of matrices $\{A_I\}=\{\mathcal A_i\}\cup \{ {\mathbb I}_{n\times n}\}$. If for any  $(m+1)$-tuple of vectors $v_I$ the corresponding quadratic map 
\bea
f_I=x^* A_{I\,} x-v_I^* x-x^* v_I\ ,
\eea
satisfies the convexity criterion of Proposition 2, $z_{\rm max}\rightarrow \infty$, then ${\mathcal F}(\mathcal A)$ is convex.

The  two new sufficient conditions for convexity formulated above invoke auxiliary map $f$. We leave the task of formulating these criteria for convexity of ${\mathcal F}(\mathcal A)$ free of any reference to $f$ and $A_i,v_i$ for the future.

\section{Different Approaches to Calculating  $z_{\rm max}$} 
\label{zmax}
In section \ref{convexityFf} we proved that a compact part of ${\mathscr F}(f)$ lying in a half-space defined by the hyperplane  ${\mathscr H}_{c_+}(z_0+z_{\rm max})$ is compact. To make this a practical method of carving a compact subregion within ${\mathscr F}(f)$ suitable for applications we would need to be able to determine $z_{\rm max}$ for a given $c_+$. A straightforward approach would be to use the definition \eqref{MR4}. We prefer to rewrite $c\cdot A -\lambda^+_{\rm min}(c\cdot A)A_+$ as ${\bf c}\cdot A$, where ${\bf c}(c)$ is given by \eqref{cdef}. Notice that matrices ${\bf c}\cdot A$ for all $c\in \R^m,\ c\cdot c_+=0$, exhaust all positive-semidefinite  combinations of $A_i$ with a non-trivial kernel. Hence vectors $c$ orthogonal to $c_+$ is enough to parametrize the entire boundary of the convex cone ${\mathcal K}^+$ of the positive-semidefinite linear combinations of $A_i$. Furthermore $c\cdot (v-A x_0)={\bf c}\cdot (v-A x_0)$ and therefore minimization problem \eqref{MR4}  can be naturally defined on the boundary $\partial {\mathcal K}^+$,
\bea
\label{dK}
z_{\rm max}=\lim\limits_{\epsilon\rightarrow 0^+} \min\limits_{c\in \partial {\mathcal K}^+} |(c\cdot A+\epsilon)^{-1}c\cdot (v-A_{}x_0)|_+^2\ ,\\
{\mathcal K}^+=\{c: c\in\R^m,\ c\cdot A\succeq 0 \}\ .
\eea 
Unfortunately this problem is not convex. A naive extension of this problem to the interior of ${\mathcal K}^+$ is not a viable option because  $\lim\limits_{\epsilon\rightarrow 0^+} \min\limits_{c\in {\mathcal K}^+} |(c\cdot A+\epsilon)^{-1}c\cdot (v-A_{}x_0)|_+^2=0$,  even when $z_{\rm max}>0$. We conclude that the minimization problem \eqref{MR4} can not be immediately reduced to a convex optimization problem admitting an efficient solution. 

This motivates us to approach the problem of calculating $z_{\rm max}$ from a slightly different angle. As was discussed in section~\ref{geomeaning} the expression minimized in \eqref{dK} is finite in the limit $\epsilon\rightarrow 0$ only when for some vector $c\neq 0$ such that $c\cdot A\succeq 0$ and $c\cdot A\not\succ 0$, i.e.~$c\cdot A$ is degenerate, the following equation has a solution, $(c\cdot A_{\,})x=c\cdot (v-Ax_0)$. The latter is trivially equivalent to solvability of $(c\cdot A_{\,})x=c\cdot v$.
Let us define the set of all such vectors $c$,
\bea
\label{Cminus}
C_-=\{c: c\in\R^m,\ |c|^2=1,\ c\cdot A\succeq 0,\ c\cdot A\not\succ 0, \ \exists\, x,\ (c\cdot A)x=c\cdot v \}.\ 
\eea
Because of the additional constraints, in real case  $C_-\subset \partial {\mathcal K}^+$ is a $(m-3)$-dimensional  sub-manifold of $(m-1)$-dimensional $\partial {\mathcal K}^+$. (In complex case  $C_-$ is $(m-4)$-dimensional.)  In principle $z_{\rm max}$ can be calculated by minimizing 
\bea
\label{zc}
z(c)= |(c\cdot A)^{-1}c\cdot v -x_0|_+^2\ ,
\eea 
over $C_-$,
\bea
\label{zmaxN}
z_{\rm max}=\min_{c\in C_-} z(c)\ .
\eea
In \eqref{zc} $(c\cdot A)^{-1}$ denotes pseudo-inverse. This leads to the following application -friendly reformulation of Proposition 2.\\ 
\noindent{\bf Proposition 2'.} The compact subset of the image of the quadratic map \eqref{fullimage},
\bea
{\mathscr F}(f,c_+,z_{\rm max})=\{y: y\in {\mathscr F}(f),\ z_0\le c_+\cdot y\le z_0+z_{\rm max} \}\subset {\mathscr F}(f),
\eea
with $z_{\rm max}$ defined in (\ref{zmaxN}, \ref{zc}), 
 is convex. (See \eqref{cplus} for the definition of $z_0$.)  If $z_{\rm max}$ is infinite the whole image ${\mathscr F}(f)$ is convex. \\
\noindent {\it Comment 4.} The sufficient condition for convexity can be concisely formulated as follows: if the set $C_-$ is empty, the image ${\mathscr F}(f)$ is convex.
  When the matrix $c\cdot A$ is degenerate, the solvability of $(c\cdot A)x=c\cdot v$ is equivalent to $x_{\rm null}^* (c\cdot v)=0$ for all $x_{\rm null}$ satisfying $(c\cdot A)x_{\rm null}=0$. Hence, if (i) a quadratic map \eqref{QM} is definite, i.e.~there is $c_+$ such that $A_+:=c_{+}\cdot A\succ 0$ and (ii) for any vector $c\neq 0$, such that $c\cdot A$ is positive semi-definite and degenerate, there is a vector $x_{\rm null}$ satisfying $(c\cdot A)x_{\rm null}=0$ and $x_{\rm null}^* (c\cdot v)\neq 0$, then the image of $f(x)$ is convex.  

When the set $C_-$ is non-empty and the map is definite, the image  ${\mathscr F}(f)$ is not strictly convex and in a general case it is non-convex at all. In this case one potential strategy to calculate  $z_{\rm max}$ would be to minimize \eqref{zmaxN} numerically, e.g.~via gradient descent along $C_-$.

\subsection{Conservative Estimate of $z_{\rm max}$}
\label{aprxsection}
Calculating $z_{\rm max}$ exactly could be a difficult. Nevertheless for many practical application it would be enough to have an easy-to-calculate conservative estimate $z_{\rm est}\le z_{\rm max}$. A very similar problem of estimating $\varepsilon_{\rm max}^2$ \eqref{MY} was addressed in \cite{my} and here we employ the same strategy.  The first step is the inequality 
\bea
\label{aprx}
z_{\rm est}:= \min_{c\in {\mathcal C}}{|c\cdot \tilde v|^2\over ||c\cdot \tilde A-\lambda_{\rm min}(c\cdot \tilde A) ||^2}\le z_{\rm max}\ , \\
{\mathcal C}=\{c: c\in \R^m,\ |c|^2=1,\ c\cdot c_+=0\}\ ,
\eea
where $\tilde v_{i}=\mathcal O (v_i-A_{i\,}x_0)$, $\tilde{A}_i=\mathcal O_{\,} A_{i\,} \mathcal O^*$ and $\mathcal O$ is defined through $A_+=(\mathcal O^*\mathcal O)^{-1}$.
Because \eqref{aprx} is homogeneous in $c$ the condition $|c|^2=1$ can be substituted by $c^* g_{} c=1$ for some $g$ which is  positive-definite on the orthogonal compliment to $c_+$ within $\R^m$. Let's choose $g_{ij}=\Re\left(\tilde{v}_i^* \tilde v_j\right)$ which satisfies this requirement (obviously $g_{ij}$ is non-negative; if it develops a zero eigenvalue on the orthogonal compliment to $c_+$, then, alas, $z_{\rm max}$ would be zero anyway). It is convenient to diagonalize $g$ and bring it to the standard form $\Lambda g_{} \Lambda^T={\rm diag}(1,\dots,1,0)$ with the help of some non-degenerate real-valued $m\times m$ matrix  $\Lambda$. Now we can define $\hat{A}_i=\sum_j \Lambda_i^j \hat A_j$ and since $\Lambda^i_{m}\sim c_{+}^i$ we introduce a new $(m-1)$-dimensional vector $\hat{c}$ with the components $\hat c=(\hat c_1,\dots,\hat c_{m-1})^T$,
\bea
\label{aprxmax}
z_{\rm est}= \left(\max_{|\hat c|^2=1}\, \left( \lambda_{\rm max}(\hat c\cdot \hat{A})-\lambda_{\rm min}(\hat c\cdot \hat{A})\right)^2\right)^{-1}\ ,\quad \hat c\in\R^{m-1}\ .
\eea
Furthermore using the inequality $\lambda_{\rm max}(A)-\lambda_{\rm min}(A)\le 2\max\{\lambda_{\rm max}(A),\lambda_{\rm max}(-A)\}$ we arrive at the following approximate conservative estimate of $z_{\rm est}$,
\bea
{\rm z}=\left(\max_{|\hat c|^2=1}\,  2_{}\lambda_{\rm max}(\hat c\cdot \hat{A})\right)^{-2} \le z_{\rm est}\ . \label{silly}
\eea  
Next step would be to calculate or estimate the Lipschitz constant $\max_{|\hat c|^2=1} \lambda_{\rm max}(\hat c\cdot \hat A)$, the problem which was previously  addressed in \cite{my}. 

In fact the estimate \eqref{silly} can be improved if we notice that $\lambda_{\rm max}(A)-\lambda_{\rm min}(A)$ is invariant under the shifts of $A$ by the identity matrix $A\rightarrow A+\mu_{\,} {\mathbb I}_{n\times n}$. Hence the better estimate for $z_{\rm est}$ would be 
\bea
{\rm z}=4^{-1}\left(\min_{\mu}\max_{|\hat c|^2=1}\,  \left(\lambda_{\rm max}(\hat c\cdot \hat{A})+\mu(\hat c)\right)\right)^{-2} \le z_{\rm est}\ ,
\eea 
where the minimum is taken over the space of functions $\mu(\hat{c})$ satisfying $\mu(-\hat{c})=-\mu(\hat{c})$. Obviously we can restrict the class of functions $\mu$ by paying the price of somewhat deteriorating the quality of the estimate. For example $\mu$ could be chosen to be a linear function $\mu(\hat{c})=\hat{c}\cdot \mu$ defined by a  vector $\mu_i$. Depending on the chosen method to estimate $\max_{|\hat c|^2=1} \lambda_{\rm max}(\hat c\cdot (\hat A+\mu))$ we can either find a minimum with respect to $\mu_i$ analytically or leave it to numerical analysis.  

In case minimization with respect to $\mu_i$ is difficult, one guideline to choose $\mu_i$ could be the following. 
The inequality $\lambda_{\rm max}(A)-\lambda_{\rm min}(A)\le 2\max\{\lambda_{\rm max}(A),\lambda_{\rm max}(-A)\}$ is saturated when $\lambda_{\rm max}(A)=-\lambda_{\rm min}(A)$ and therefore it makes sense to choose $\mu_i$ such that any combination $\hat c\cdot (\hat A+\mu)$ has both positive and negative eigenvalues. For example this can be done with help of $\mu_i=-\Tr(\hat{A}_i)$, which  is an optimal condition for $n=2$.
%

\section{Examples and Applications} 
\label{examples}
\subsection{Convexity of solvability set of Power Flow equations for DC networks}
\label{proofPF}
The Power Flow equations for DC networks is a set of quadratic equations which express power injections on each bus of the network in terms of real-valued voltages. For a system consisting of $N$ buses, $i=1,\dots,N$, the equations take the following form 
\bea
\label{PF}
P_i=V_i \left(\sum_{j=1}^N Y_{ij} V_j\right)\ ,\quad i=1,\dots,N\ ,
\eea
where $Y_{ij}$ is a symmetric $N\times N$ admittance matrix of resistances which is defined for each network as follows 
\begin{eqnarray}
\label{Y}
Y_{ik}&=&\left\{
\begin{array}{c r}
\sum_{l\sim i} y_{il} & {\rm \ if\ }i=k\\
-y_{ik} & {\rm \ if\ }i\sim k\\
0 & {\rm \ if\ }i\not\sim k\\
\end{array}
\right.
\end{eqnarray}
Here we adopt the notation $i\sim k$ to indicated that nodes (buses) $i$ and $k$ are connected by a line, and in \eqref{Y} we assume that each line is purely resistive, $y_{ik}=y_{ki}>0$, for any $i\sim k$. In what follows we will also assume the network is connected i.e.~any two nodes can be connected by a combination of lines. 

We choose  $i=N$ to be the slack bus, which means $V_N\equiv 1$. Then 
the Power Flow equations should be understood as equations  expressing powers $P_i$ for $i=1,\dots,N-1$, in terms of voltages $V_i$, where  $i=1,\dots,N-1$. Hence the Power Flow equations for DC networks define a ${\mathbb R}^{N-1}\rightarrow {\mathbb R}^{N-1}$ quadratic map. The solvability set of Power Flow equations is the combination of all feasible powers $P_i$ such that equations \eqref{PF} admit a solution, i.e.~it is the image of this quadratic map. 

In what follows we prove that solvability set of Power Flow equations \eqref{PF} is convex. This result is closely related to the proof that convex relaxation for Optimal Power Flow problem for DC networks is exact \cite{Low}. At the same time  convexity of solvability set is a new result which was not covered in literature previously.

The quadratic map defined by \eqref{PF} can be represented in the conventional form \eqref{QM} with $n=m=N-1$, $f_i^0=0$, while $A$, $v$ and $x$ are given by
\bea
\label{PFmap}
(A_i)_{kl}={1\over 2}\left(\delta_{ik}Y_{il}+\delta_{il}Y_{ik}\right)\ , \quad
v_i=-{1\over 2}Y_{Ni}\ ,\quad x^i=V_i\ ,\quad  i,k,l=1,\dots,n.
\eea
To verify that the set of matrices $A_i$ is definite it is enough to calculate the sum
\bea
\sum_{i=1}^n x^* A_i x =\left.\sum_{i=1}^N P_i\right|_{V_N=0}={1\over 2}\sum_{i\sim j} y_{ij}(V_i-V_j)^2+\sum_{i=1}^n y_{Ni}V_i^2\ ,
\label{Pplus}
\eea
which is manifestly positive-definite. From the physics point of view, definiteness of map \eqref{PFmap} is a consequence of the network being resistive, hence the total power dissipated by the lines (first term in \eqref{Pplus}) is non-negative. From here we immediately conclude that vector $c_+$ in \eqref{cplus} can be chosen to be 
\bea
c_+=(1,\dots,1)^T\ .
\eea
Next we would like to identify the set $C_-$ \eqref{Cminus}.  We assume that vector $c\in C_-$ which in particular means that matrix $(c\cdot A)$ is positive semi-definite and degenerate,
\bea
(c\cdot A)_{ik}=\left\{
\begin{array}{c r}
c_i\sum_{l\sim i} y_{il} & {\rm \ if\ }i=k\\
-(c_i+c_k)y_{ik}/2 & {\rm \ if\ }i\sim k\\
0 & {\rm \ if\ }i\not\sim k\\
\end{array}
\right.
\eea
Since $(c\cdot A)\succeq 0$ all diagonal elements of $(c\cdot A)$ have to be non-negative and since all $y_{ik}>0$ we find that all $c_i\ge 0$. In fact all $c_i$ must be strictly positive. Indeed, let us assume that  $c_i=0$ for a particular $i$. For $c\cdot A$ to be positive semi-definite all elements $(c\cdot A)_{ik}$  must be zero, $(c_i+c_k)y_{ik}=0$. This is only possible if $c_k=0$ for all $k\sim i$. And since the network is connected, repeating this consideration enough times, we find $c_k=0$ for all $k=1,\dots,n$. This contradicts the assumption $|c|^2=1$, and therefore all $c_i>0$.

By definition of $C_-$ matrix $c\cdot A$ must be singular, i.e.~there is a normalized 
 vector $x_{\rm null}$ such that $(c\cdot A)x_{\rm null}=0$. This vector minimizes the  following quadratic form 
\bea
\label{QF}
x_{\rm null}^*(c\cdot A)x_{\rm null}= \sum_{i=1}^n(c\cdot A)_{ii} V_i^2 -2\sum_{i\sim k} (c_i+c_j)y_{ik} V_i V_k\ ,
\eea
subject to $|x_{\rm null}|^2=\sum_i V_i^2=1$. Taking into account that for any $i\sim k$ the combination $ (c_i+c_j)y_{ik}>0$, for the given values of $|V_i|$, 
the quadratic form (\ref{QF}) will be minimal if $V_i V_k=|V_i||V_k|$. Hence all components of $x_{\rm null}$ minimizing (\ref{QF}) must have the same sign.  Without loss of generality we can assume that $V_i\ge 0$. In fact it is possible to show that all $V_i$ minimizing \eqref{QF} must be strictly positive (or negative). To prove that, we assume the opposite, that $V_i=0$ for some $i$. Then we analyze the equation $(c\cdot A)x_{\rm null}=0$, namely look at the $i$-the component of vector $(c\cdot A)x_{\rm null}$,
\bea
c_i\sum_{k\sim i}y_{ik} V_i-\sum_{k\sim i}(c_i+c_k)y_{ik} V_k/2=0\ .
\eea
Since all $V_k$ are non-negative and by assumption $V_i=0$ this equation 
can be satisfied  only if $V_k=0$ at all nodes $k$ connected with $i$. Since the network is connected, by repeating this logic, we find that all components of vector $x_{\rm null}$ minimizing \eqref{QF} are zero, which contradicts the assumption $|x_{\rm null}|^2=1$. 

Finally, we consider if the equation  $(c\cdot A)x=c\cdot v$ can have a solution. The necessary and sufficient condition for the solution to exist is 
\bea
\label{condition2}
x^*_{\rm null}(c\cdot v)=0
\eea
for all vectors $x_{\rm null}$ satisfying $(c\cdot A)x_{\rm null}=0$.
The equation \eqref{condition2} written in terms of $V_i$ has the form 
\bea
\sum_{i\sim N} c_i y_{Ni} V_i=0\ .
\eea
Clearly it can not be satisfied because it is a sum of strictly positive terms. We therefore conclude that $C_-$ is empty, hence $z_{\rm max}$ is infinite and according to Proposition 2' the solvability set is convex. 

\subsection{3-bus example}
\label{3bus}
To illustrate some of the ideas discussed above, let us consider a simple network consisting of $3$ buses all connected to each other. The corresponding admittance matrix and the quadratic map is given by 
\bea
Y&=&\left(
\begin{array}{ccc}
2 & -1/2 & -3/2 \\
-1/2 & 2 & -3/2\\
-3/2 & -3/2 & 3
\end{array}
\right)\ ,\\
\label{P1}
P_1&=&V_1(2V_1-V_2/2-3/2V_3)\ ,\\
P_2&=&V_2(2V_2-V_1/2-3/2V_3)\ .
\label{P2}
\eea
We have established in section \ref{proofPF} that the image of the quadratic map $f_i=P_i(V_1,V_2)$, $i=1,2$, with $V_3$ being a constant is convex. 
In section \ref{connection} we establish that the convexity of ${\mathscr F}(f)$ implies convexity of image ${\mathscr F}({\mathsf f})$ of the homogeneous quadratic map ${\mathsf f}$ defined through \eqref{FFH}. This simply means the functions (\ref{P1},\ref{P2}) are amended by 
\bea
P_3=V_3^2\ ,
\eea
and ${\mathsf f}_I=P_I(V_1,V_2,V_3)$ are the functions of all three variables $V_i$. The associated matrices ${\mathsf A}_I$ are as follows
\bea
{\mathsf A}_1=\left(
\begin{array}{ccc}
 2 & -\frac{1}{4} & -\frac{3}{4} \\
 -\frac{1}{4} & 0 & 0 \\
 -\frac{3}{4} & 0 & 0 \\
\end{array}
\right),\ \
{\mathsf A}_2=\left(
\begin{array}{ccc}
 0 & -\frac{1}{4} & 0 \\
 -\frac{1}{4} & 2 & -\frac{3}{4} \\
 0 & -\frac{3}{4} & 0 \\
\end{array}
\right),\ \
{\mathsf A}_3=\left(
\begin{array}{ccc}
 0 & 0 & 0 \\
 0 & 0 & 0 \\
 0 & 0 & 1 \\
\end{array}
\right).
\eea
The sum ${\mathsf A}_+={\mathsf A}_1+{\mathsf A}_2+{\mathsf A}_3$ is positive-definite. To find matrices ${\mathcal A}$ and the joint numerical range ${\mathcal F}({\mathcal A})$ associated with the base of the cone ${\mathscr F}({\mathsf f})$ we consider two linear combinations of ${\mathsf A}_I$ which are linearly independent with ${\mathsf A}_+$,
\bea
{\mathcal A}_1={{\mathsf A}_2-{\mathsf A}_1\over \sqrt{2}}\ ,\quad {\mathcal A}_2={{\mathsf A}_3\over 4}\ .
\label{Amatr}
\eea
Upon changing the coordinates in the space of $V_i$ we can bring ${\mathsf A}_+$ to be the identity matrix. In these coordinates matrices ${\mathcal A}_i$, $i=1,2$, \eqref{Amatr} are given in the caption of Fig.~\ref{fig:flatedge}. The convexity of the corresponding joint numerical range depicted in Fig.~\ref{fig:flatedge} is thus guaranteed by the {\it first} sufficient condition for convexity formulated in section \ref{newc}.  It should be noted for completeness that  since ${\mathcal F}({\mathcal A})$ in this case is only $2$-dimensional, its convexity is also guaranteed by the result of \cite{Brickman}.

\subsection{4-bus example}
\label{4bus}
Our next example is a $4$-bus network with the admittance matrix 
\bea
\label{YY}
Y=\left(
\begin{array}{cccc}
 0 & 1 & 0 & 1 \\
 1 & 0 & 1 & 1 \\
 0 & 1 & 0 & 1 \\
 1 & 1 & 1 & 0 \\
\end{array}
\right)\ .
\eea
The corresponding expressions $f_i=P_i(V_1,V_2,V_3,V_4=1)$ \eqref{PF} define the initial quadratic map $f$. The latter gives rise to the homogeneous quadratic map ${\mathsf f}$ \eqref{FFH} with matrices  ${\mathsf A}_I$, $I=1,\dots,4$, given by 
\bea
\nonumber
{\mathsf A}_1&=\left(
\begin{array}{cccc}
 2 & -\frac{1}{2} & 0 & -\frac{1}{2} \\
 -\frac{1}{2} & 0 & 0 & 0 \\
 0 & 0 & 0 & 0 \\
 -\frac{1}{2} & 0 & 0 & 0 \\
\end{array}
\right), \quad 
{\mathsf A}_2&=\left(
\begin{array}{cccc}
 0 & -\frac{1}{2} & 0 & 0 \\
 -\frac{1}{2} & 3 & -\frac{1}{2} & -\frac{1}{2} \\
 0 & -\frac{1}{2} & 0 & 0 \\
 0 & -\frac{1}{2} & 0 & 0 \\
\end{array}
\right), \\ 
{\mathsf A}_3&=\left(
\begin{array}{cccc}
 0 & 0 & 0 & 0 \\
 0 & 0 & -\frac{1}{2} & 0 \\
 0 & -\frac{1}{2} & 2 & -\frac{1}{2} \\
 0 & 0 & -\frac{1}{2} & 0 \\
\end{array}
\right),\quad
{\mathsf A}_4&=\left(
\begin{array}{cccc}
 0 & 0 & 0 & 0 \\
 0 & 0 & 0 & 0 \\
 0 & 0 & 0 & 0 \\
 0 & 0 & 0 & 1 \\
\end{array}
\right).
\label{A4}
\eea
We choose ${\mathsf c}_+=2(1,1,1,1)^T$, one can check ${\mathsf A}_+\equiv {\mathsf c}_+\cdot {\mathsf A}\succ 0$,
and define ${\mathcal A}_i$, $i=1,2,3$, through the following linear combinations 
\bea
{\mathcal A}_1={{\mathsf A}_1-{\mathsf A}_3\over 3}\ ,\quad  
{\mathcal A}_2={{\mathsf A}_2-{\mathsf A}_3\over 3}\ ,\quad  
{\mathcal A}_3={\mathsf A}_+-{\mathsf A}_4\ .
\eea 
The joint numerical range ${\mathcal F}(\mathcal A)$ is thus the base of the cone ${\mathscr F}({\mathsf f})$. 
In the coordinates where matrix ${\mathsf A}_+$ is the identity matrix, matrices ${\mathcal A}_i$ have the following form
\bea
\nonumber
{\mathcal A}_1=\left(
\begin{array}{cccc}
 \frac{1}{4 \sqrt{2}} & 0 & 0 & \frac{1}{12 \sqrt{2}} \\
 0 & -\frac{1}{4 \sqrt{2}} & 0 & -\frac{1}{12 \sqrt{2}} \\
 0 & 0 & 0 & 0 \\
 \frac{1}{12 \sqrt{2}} & -\frac{1}{12 \sqrt{2}} & 0 & 0 \\
\end{array}
\right),\ 
{\mathcal A}_2=\left(
\begin{array}{cccc}
 \frac{1}{8 \sqrt{2}}-\frac{7}{72} & -\frac{1}{72} & \frac{1}{144} & \frac{1}{48} \left(\sqrt{2}-3\right) \\
 -\frac{1}{72} & -\frac{7}{72}-\frac{1}{8 \sqrt{2}} & \frac{1}{144} & \frac{1}{48} \left(-\sqrt{2}-3\right) \\
 \frac{1}{144} & \frac{1}{144} & \frac{7}{36} & \frac{1}{8} \\
 \frac{1}{48} \left(\sqrt{2}-3\right) & \frac{1}{48} \left(-\sqrt{2}-3\right) & \frac{1}{8} & 0 \\
\end{array}
\right),\ \ 
\eea
\bea
\label{A3}
{\mathcal A}_3&=&\left(
\begin{array}{cccc}
 1 & 0 & 0 & 0 \\
 0 & 1 & 0 & 0 \\
 0 & 0 & 1 & 0 \\
 0 & 0 & 0 & -1 \\
\end{array}
\right)\ .
\eea
According to the {\it first} sufficient condition for convexity formulated in section \ref{newc}, the joint numerical range  ${\mathcal F}(\mathcal A)\subset{\mathbb R}^3$  is convex. We would like to emphasize that this is a new result as the convexity of such joint numerical range could not have been established using previously known criteria. 
\begin{figure}[t]
\includegraphics[width=0.32\textwidth]{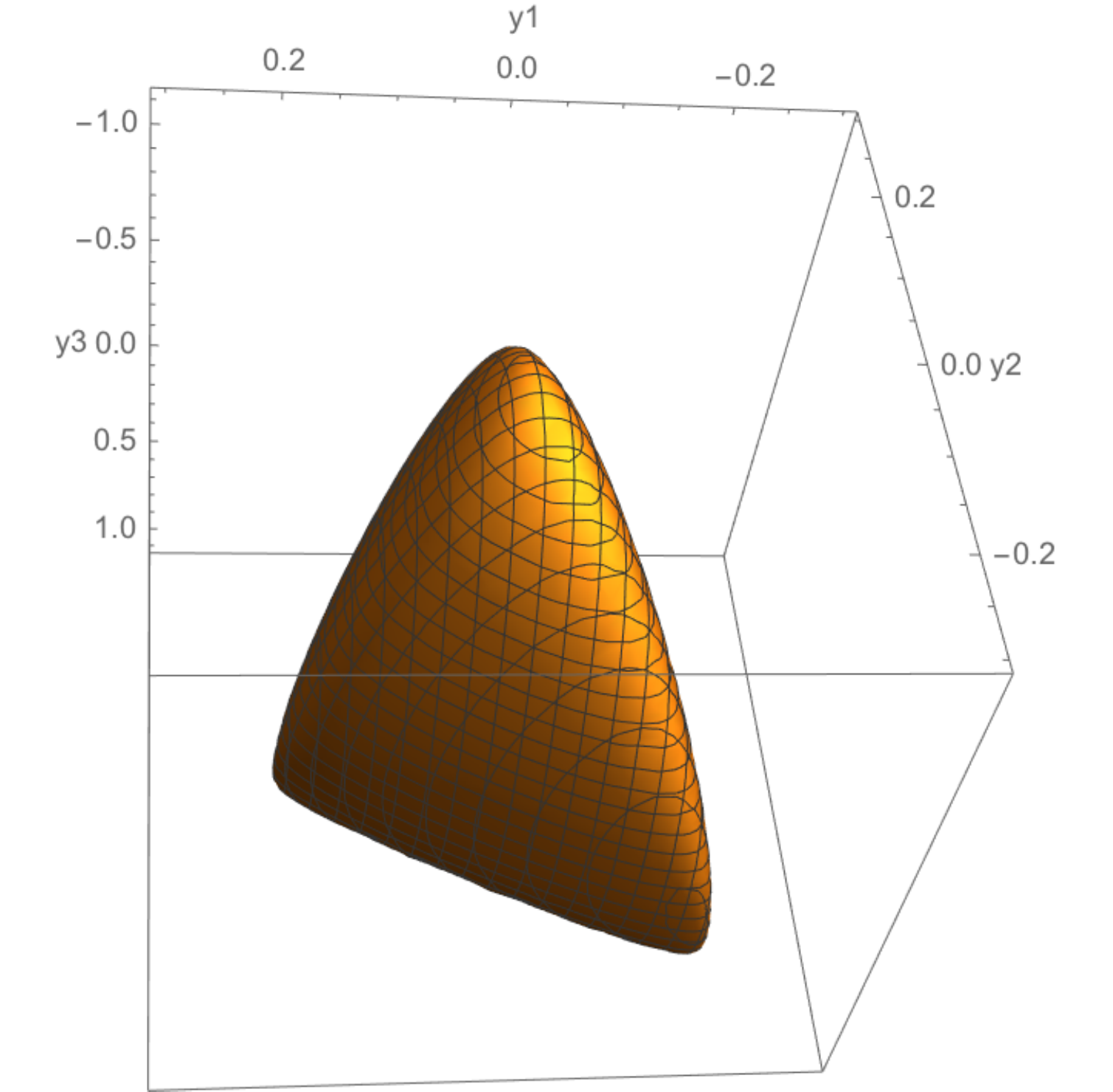}
\includegraphics[width=0.32\textwidth]{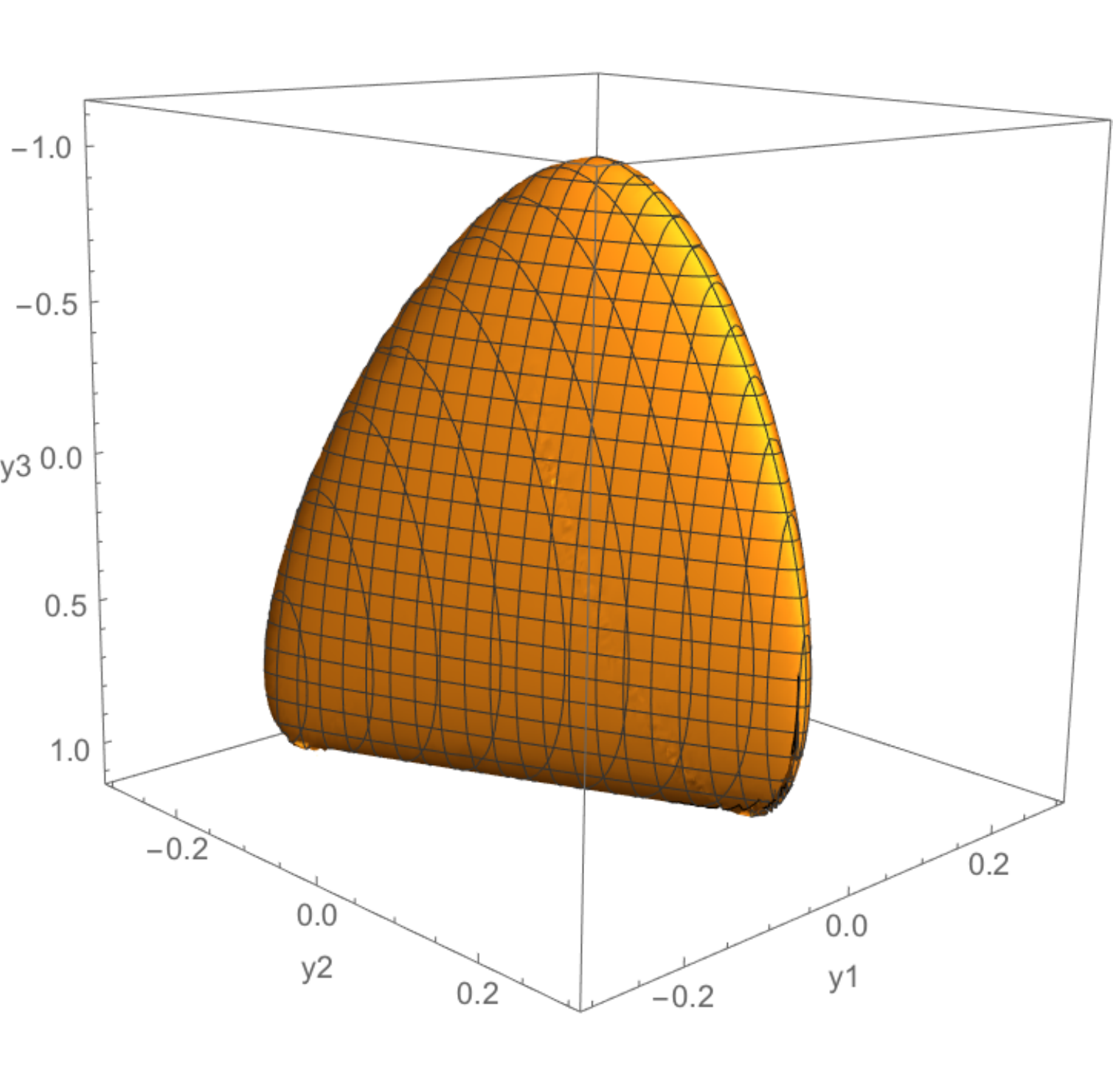}
\includegraphics[width=0.32\textwidth]{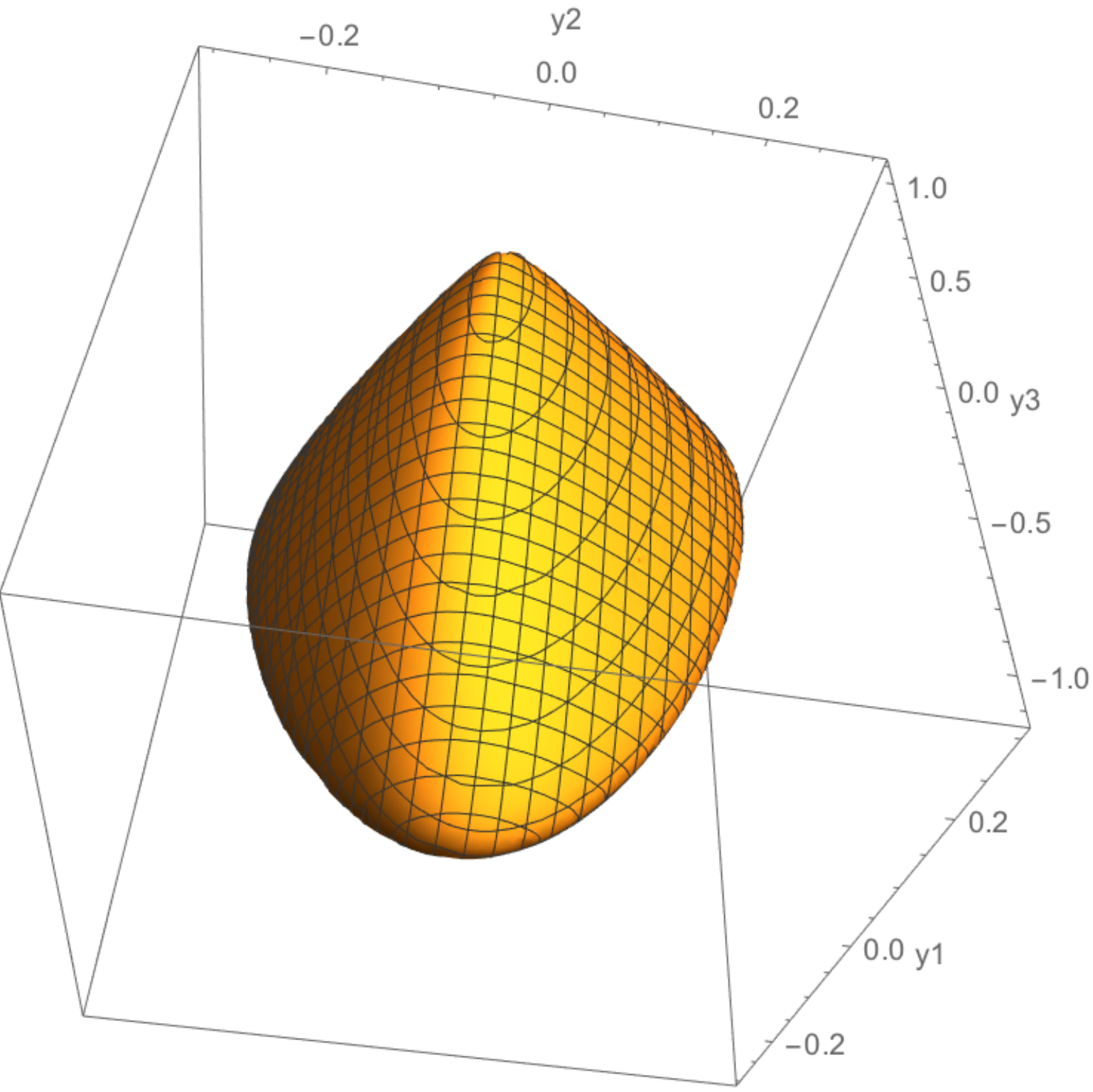}
\caption{Different projections of the joint numerical range ${\mathcal F}({\mathcal A})$ of the matrices \eqref{A3}. This joint numerical range satisfies the new criterion for convexity formulated in section \ref{newc}.}
\label{3Dconvex}
\end{figure}

\begin{SCfigure}[50][t]
\includegraphics[width=0.4\textwidth]{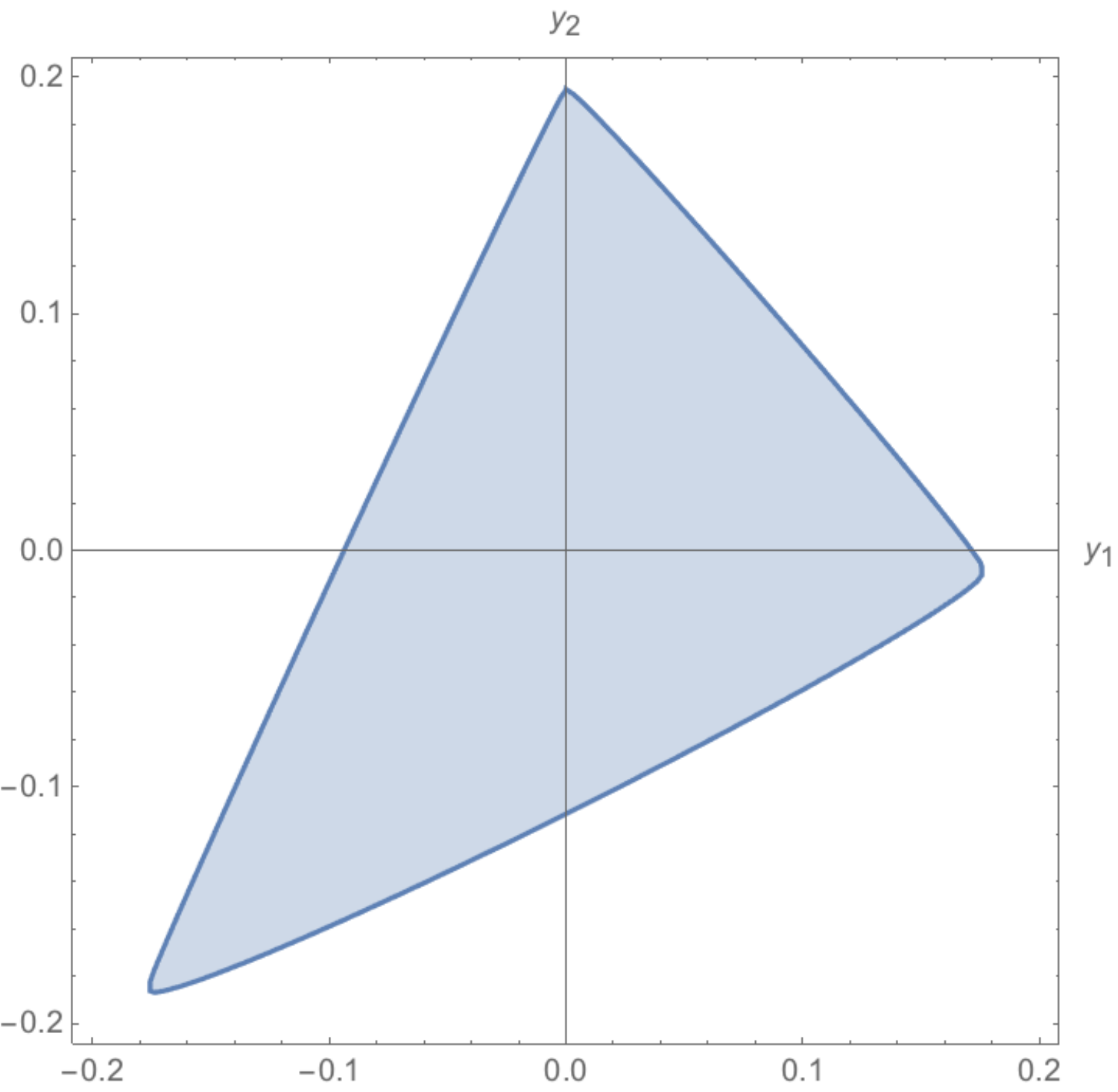}
\captionsetup[top]{singlelinecheck=on}
\caption[top]{The joint numerical range ${\mathcal F}({\mathcal A})\subset {\mathbb R}^2$, i.e.~the image of $y_i(x)=x^* {\mathcal A}_i x$, $|x|^2=1$, with ${\mathcal A}_i$, given by \eqref{a3}. It is the ``bottom" of the joint numerical range ${\mathcal F}({\mathcal A})\subset{\mathbb R}^4$ with ${\mathcal A}_i$ given by \eqref{A3}, namely it is an intersection of ${\mathcal F}({\mathcal A})$ with the supporting hyperplane $y_3=1$.\vspace{2cm}}
\label{bottom}
\end{SCfigure}

The corresponding joint numerical range $y_i(x)=x^* \mathcal A_i x$, $i=1,2,3$, $|x|^2=1,\ x\in {\mathbb R}^3$, is shown in Fig.~\ref{3Dconvex}. The ``bottom" of that figure, namely the intersection of the image of $y(x)$ with the hyperplane $y_3=1$, is a flat two-dimensional convex set which is shown in Fig.~\ref{bottom}. 

Convexity of the latter can be independently deduced from the  {\it second} new sufficient condition for convexity formulated in section \ref{newc}. Namely, one can start with the  expression for powers \eqref{PF} specified by the admittance matrix \eqref{YY}. Assuming $V_4=0$, powers $P_i$, $i=1,2,3$, become quadratic homogeneous functions of $V_i$,  $i=1,2,3$. The corresponding matrices $A_i$ are given by the $3\times 3$ upper left sub-matrices of the matrices \eqref{A4},
\bea
A_1=\left(
\begin{array}{ccc}
 2 & -\frac{1}{2} & 0  \\
 -\frac{1}{2} & 0 & 0  \\
 0 & 0 & 0 
\end{array}
\right), \ \
A_2&=\left(
\begin{array}{cccc}
 0 & -\frac{1}{2} & 0  \\
 -\frac{1}{2} & 3 & -\frac{1}{2} \\
 0 & -\frac{1}{2} & 0 
\end{array}
\right),\ \  
A_3&=\left(
\begin{array}{cccc}
 0 & 0 & 0  \\
 0 & 0 & -\frac{1}{2}  \\
 0 & -\frac{1}{2} & 2  
\end{array}
\right) ,\nonumber
\eea
while $v_i=0$ and $f_i^0=0$. This homogeneous quadratic map is definite, $A_+=2(A_1+A_2+A_3)\succ 0$. By performing a linear change of variables $V_i$ we can bring $A_+$ to be the identity matrix, while two linearly independent combinations ${\mathcal A}_1=(A_1-A_3)/3$ and ${\mathcal A}_2=(A_2-A_3)/3$ will be given by (compare with \eqref{A3}),
\bea
{\mathcal A}_1=\left(
\begin{array}{ccc}
 \frac{1}{4 \sqrt{2}} & 0 & 0 \\
 0 & -\frac{1}{4 \sqrt{2}} & 0 \\
 0 & 0 & 0 \\
\end{array}
\right),\ \ 
{\mathcal A}_1=\left(
\begin{array}{ccc}
 -\frac{7}{72}+\frac{1}{8 \sqrt{2}} & -\frac{1}{72} & \frac{1}{144} \\
 -\frac{1}{72} & -\frac{7}{72}-\frac{1}{8 \sqrt{2}} & \frac{1}{144} \\
 \frac{1}{144} & \frac{1}{144} & \frac{7}{36} \\
\end{array}
\right) .
\label{a3}
\eea
The convex image of $y_i=x^* {\mathcal A}_i x$, $i=1,2$, $|x|^2=1,\ x\in {\mathbb R}^3$ is shown in Fig.~\ref{bottom}.  

\subsection{${\mathbb C}^2\rightarrow  {\mathbb R}^4$ example}
Consider a quadratic map $f:{\mathbb C}^2\rightarrow  {\mathbb R}^4$ with $A_i$ being Pauli matrices
\bea
&A_1=\left(
\begin{array}{cc}
 0 & 1 \\
 1 & 0 \\
\end{array}
\right),\quad 
&A_2=\left(
\begin{array}{cc}
 0 & -i \\
 i & 0 \\
\end{array}
\right) ,\nonumber \\
&A_3=\left(
\begin{array}{cc}
 1 & 0 \\
 0 & -1 \\
\end{array}
\right) ,\quad 
&A_4=\left(
\begin{array}{cc}
 1 & 0 \\
 0 & 1 \\
\end{array}
\right)\ ,
\label{Pa}
\eea
and $v$ given by 
\bea
\label{Pv}
v_1=\left(\begin{array}{c} 1\\ 0 \end{array}\right),\ \ v_2=\left(\begin{array}{c} -i \\  0 \end{array}\right), \ \ v_3=\left(\begin{array}{c} 0 \\  0 \end{array}\right), \ \ 
v_4=\left(\begin{array}{c} 2 \\  0 \end{array}\right).
\eea
Clearly, this map is definite. We choose $c_+=(0,0,0,1)$ which results in $x_0=v_+=v_4$ 
and $z_0=-v_+^* A_+^{-1} v_+=-4$. Now we identify all vectors $c$ such that $c\cdot A\succeq 0$ and degenerate. First, ${\rm det}(c\cdot A)=0$ requires $c_1^2+c_2^2+c_3^2=c_4^2$, then $c\cdot A\succeq 0$ also specifies $c_4>0$. Up to an overall rescaling we parametrize all such vectors by points on ${\mathbb S}^2$,
\bea
c_1=\sin\theta \cos\phi,\quad c_2=\sin\theta \sin\phi,\quad c_3=\cos\theta,\quad c_4=1\ .
\eea
A normalized vector $x_{\rm null}$ satisfying $(c\cdot A)x_{\rm null}$ is given by 
\bea
x_{\rm null}=\left(\begin{array}{c} \sin{\theta\over 2}e^{-i\phi}\\ -\cos{\theta\over 2} \end{array}\right)\ .
\eea
The equation $x_{\rm null}^*(c\cdot v)$ yields 
\bea
-\sin \frac{\theta }{2} \left(\sin (\theta )+2 e^{i \phi }\right)=0,
\eea
with the only solution $\theta=0$. Hence $C_-$ consists of only one vector $c_-=(0,0,1,1)^T/\sqrt{2}$. It is straightforward to calculate $z(c_-)=1$ using the definition \eqref{zc}. Since $C_-$ consists of only one vector $z_{\rm max}=1$. 

Let us consider the image of $x_{\rm b}+ t x_{\rm null}$
where $x_{\rm b}=(1,0)^T$ is a solution of $(c_-\cdot A)x_{\rm b}=c_-\cdot v$ and $t$ is an arbitrary complex number. Points $y(t)=f(x_{\rm b}+ t x_{\rm null})$ belong to the ``flat edge" which is an intersection of ${\mathscr F}(f)$ with the supporting hyperplane orthogonal to $c_-$. We can find these points explicitly 
\bea
\label{mapt}
y_i(t)=\left(2 (\Re(t)-1),2 \Im(t),1-|t|^2,|t|^2-3\right)^T\ .
\eea
As expected these point lie on a hyperplane $y_3+y_4=-2$. Thus \eqref{mapt} defines a map from ${\mathbb C}^2\equiv {\mathbb R}^4$ into  ${\mathbb R}^3$ with $y_3$ being a quadratic function of $y_1$ and $y_2$. Quite clearly the image of this map (a parabolic surface in ${\mathbb R}^3$ embedded into ${\mathbb R}^4$) is not convex. 

\begin{figure}[t]
\begin{subfigure}[b]{0.32\textwidth}
\includegraphics[width=1\textwidth]{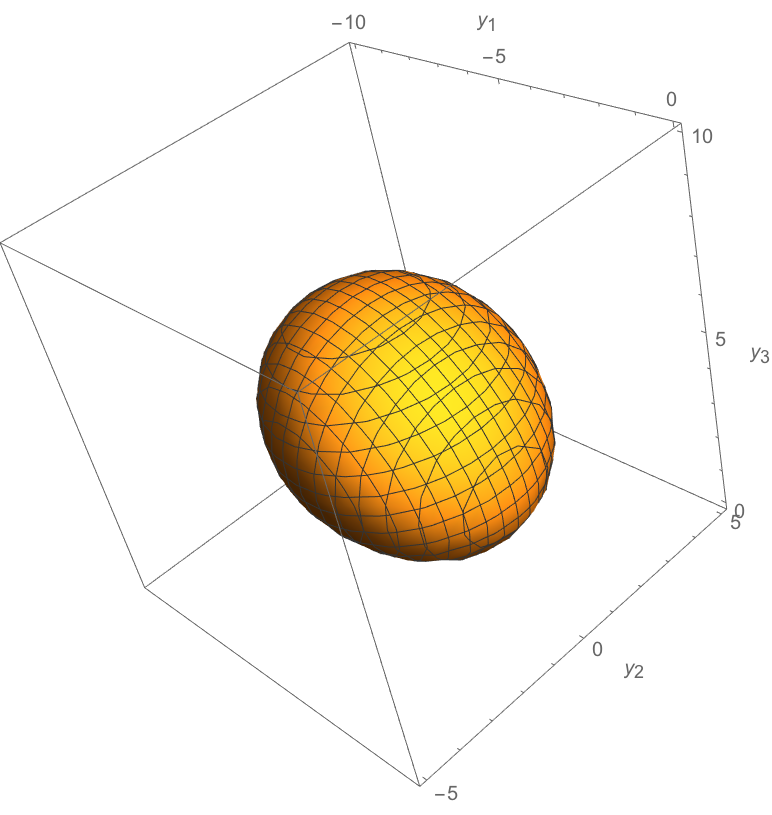}
\caption{$y_4=-3.5$}
\end{subfigure}
\begin{subfigure}[b]{0.32\textwidth}
\includegraphics[width=1\textwidth]{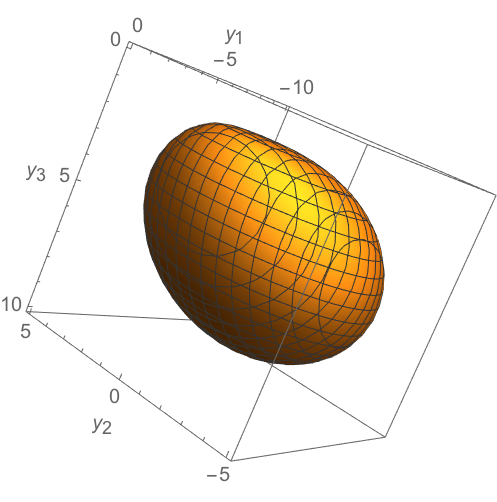}
\caption{$y_4=-3$}
\end{subfigure}
\begin{subfigure}[b]{0.32\textwidth}
\includegraphics[width=1\textwidth]{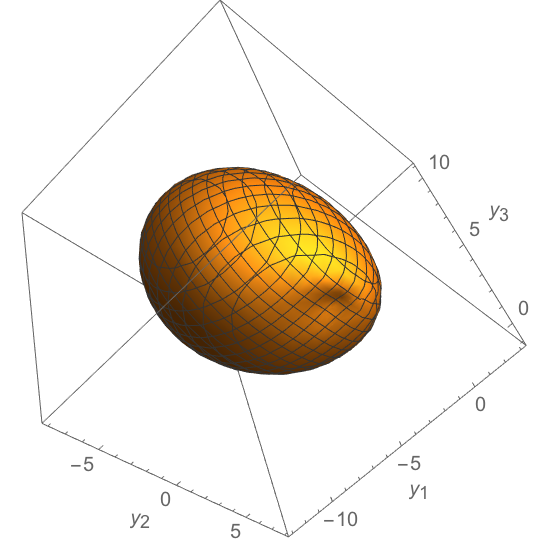}
\caption{$y_4=-2.5$}
\end{subfigure}
\caption{Intersection of the image ${\mathscr F}(f)$ of the quadratic map (\ref{Pa},\ref{Pv}) with the hyperplanes $y_4={\rm const}$. (a)  $y_4=-3.5<z_0+z_{\rm max}$ and the corresponding intersection is convex. (b)  $y_4=-3=z_0+z_{\rm max}$ and the corresponding intersection is convex but not strongly convex as it develops a ``flat edge" located at $y_3=1$. (c)  $y_4=-2.5>z_0+z_{\rm max}$ and the corresponding intersection is non-convex.}
\label{C2R4}
\end{figure}

We conclude that the image of $f(x)$ specified by (\ref{Pa},\ref{Pv}) is not convex overall, but the compact part defined by the inequality $-4=z_0\le y_4\le z_0+z_{\rm max}=-3$ is convex. This is shown in Fig.~\ref{C2R4}, where we plot the intersection of ${\mathscr F}(f)$ with the hyperplanes $y_4={\rm const}$ for different values of $y_4$.

Finally, let us estimate $z_{\rm max}$ using the approximate approach of section \ref{aprxsection}. In our case $\mathcal O$ is an identity matrix and $\tilde{v}_i=v_i-A_i x_0$ are given by 
\bea
\tilde{v}_1=\left(\begin{array}{c} 1\\ -2 \end{array}\right),\ \ \tilde{v}_2=\left(\begin{array}{c} -i \\  -2i \end{array}\right), \ \ \tilde{v}_3=\left(\begin{array}{c} -2 \\  0 \end{array}\right), \ \ 
\tilde{v}_4=\left(\begin{array}{c} 0 \\  0 \end{array}\right).
\eea
Hence matrix $g_{ij}=\Re(\tilde{v}^*_i \tilde{v}_j)$ is 
\bea
g_{ij}=\left(
\begin{array}{cccc}
 5 & 0 & -2 & 0 \\
 0 & 5 & 0 & 0 \\
 -2 & 0 & 4 & 0 \\
 0 & 0 & 0 & 0 \\
\end{array}
\right)\ .
\eea
After introducing 
\bea
\Lambda=\left(
\begin{array}{cccc}
 -\frac{\sqrt{17}+1}{\sqrt{26 \sqrt{17}+170}} & 0 & 2 \sqrt{\frac{2}{13 \sqrt{17}+85}} & 0 \\
 0 & -\frac{1}{\sqrt{5}} & 0 & 0 \\
 \frac{\sqrt{17}-1}{\sqrt{170-26 \sqrt{17}}} & 0 & 2 \sqrt{\frac{2}{85-13 \sqrt{17}}} & 0 \\
 0 & 0 & 0 & 1 \\
\end{array}
\right),
\eea 
such that $\Lambda g \Lambda^T={\rm diag}(1,1,1,0)$ we calculate  matrices $\tilde{A}_i=(\Lambda A)_i\equiv\sum_j \Lambda_i^j A_j$ for $i,j=1,2,3$. To make the presentation concise  we do not write the explicit expressions for  $\tilde{A}_i$ here.
It is important to note that since  $A_i$ for $i=1,2,3$ were traceless, so will be $\tilde{A}_i$ and $\hat{c}\cdot \hat{A}$.
 Hence $\lambda_{\rm max}(\hat{c}\cdot \hat{A})=-\lambda_{\rm min}(\hat{c}\cdot \hat{A})$ and it can be calculated explicitly 
\bea
\label{maxl}
\lambda_{\rm max}(\hat{c}\cdot \hat{A})=\frac{\sqrt{-5 \left(\sqrt{17}-9\right) \text{c1}^2+32 \text{c2}^2+5 \left(\sqrt{17}+9\right) \text{c3}^2}}{4 \sqrt{10}}\ .
\eea
Using the constraint $\sum_{i=1}^3 \tilde{c}_i^2=1$ \eqref{maxl} can be expressed as a function of two independent variables $c_1^2,c_2^2\ge 0$ with the maximum corresponding to $c_1^2=c_2^2=0$ and $c_3^2=1$,
\bea
\max_{|\tilde{c}|^2=1}\lambda_{\rm max}(\hat{c}\cdot \hat{A})={(9 + \sqrt{17})^{1/2}\over 2^{1/2}4}\ ,
\eea
and according to \eqref{silly} we find 
\bea
{\rm z}={9-\sqrt{17}\over 8}\approx 0.61\ .
\eea
Hence, our estimate assures convexity of approximately ${\rm z}/z_{\rm max}\approx 61\%$ of the maximal convex subset ${\mathscr F}(f,c_+,z_{\rm max})$ within ${\mathscr F}(f)$. 

\section*{Acknowledgments}
I would like to thank Vasily Pestun and Konstantin Turitsyn for useful discussions.


\end{document}